\documentclass[11pt,a4paper]{amsart}
\usepackage[utf8]{inputenc}
\usepackage[T1]{fontenc}
\usepackage{amsmath, amsthm, amssymb, amsfonts, mathtools, graphicx, multicol, array, enumerate, bbm, bm, comment, tikz-cd, hyperref, blindtext, subfiles}
\usepackage[dvipsnames]{xcolor}
\usepackage[margin=1in]{geometry}
\usepackage{hyperref}
\usetikzlibrary[decorations.pathmorphing]

\newtheorem{theorem}{Theorem}[section]
\newtheorem{definition}[theorem]{Definition}
\newtheorem{exmp}[theorem]{Example}
\newtheorem{lemma}[theorem]{Lemma}
\newtheorem{proposition}[theorem]{Proposition}
\newtheorem{corollary}[theorem]{Corollary}
\newtheorem{rem}[theorem]{Remark}

\newenvironment{remark}{\begin{rem}\rm}{\end{rem}\rm}

\DeclareMathOperator{\Aut}{\mathrm{Aut}}

\DeclareMathOperator{\Hom}{\mathrm{Hom}}
\DeclareMathOperator{\Ext}{\mathrm{Ext}}
\DeclareMathOperator{\ext}{\mathrm{ext}}

\DeclareMathOperator{\R}{\mathbb{R}}

\DeclareMathOperator{\C}{\mathbb{C}}

\DeclareMathOperator{\Pic}{\mathrm{Pic}}
\DeclareMathOperator{\rk}{\mathrm{rk}}

\DeclareMathOperator{\Gr}{\mathrm{Gr}}

\DeclareMathOperator{\gon}{\mathrm{gon}}
\DeclareMathOperator{\Cliff}{\mathrm{Cliff}}

\renewcommand{\P}{\mathbb{P}}
\renewcommand{\O}{\mathcal{O}}

\renewcommand{\L}{\left}
\renewcommand{\R}{\right}

\newcommand{\nocontentsline}[3]{}
\let\origcontentsline\addcontentsline
\newcommand\stoptoc{\let\addcontentsline\nocontentsline}
\newcommand\resumetoc{\let\addcontentsline\origcontentsline}

\author[Samir Geiger]{Samir Geiger}
\address{Samir Geiger, Department of Mathematics, Humboldt-Universität zu Berlin, Germany}
\email{geigersa@hu-berlin.de}

\author[Andreas Kretschmer]{Andreas Kretschmer}
\address{Andreas Kretschmer, Department of Mathematics, Humboldt-Universität zu Berlin, Germany}
\email{andreas.kretschmer@hu-berlin.de}

\title{On Green's Conjecture for Curves on Anticanonical Surfaces}

\begin{document}
\begin{abstract}
    We prove Green's conjecture for smooth curves $C$ on rational surfaces $X$ such that $-K_X|_C$ is effective on $C$. In addition, on birationally ruled surfaces $X$ of irregularity $q \geq 1$ carrying an effective anticanonical divisor, we show that Green's conjecture holds for all smooth curves $C$ on $X$ covering the base curve with degree at least $2$, provided the latter has only finitely many gonal pencils or the gonality of $C$ is special with respect to the ruling.
\end{abstract}
\maketitle

\section{Introduction}

\noindent Green's conjecture on syzygies of canonical curves predicts the vanishing statement
$$
K_{p,2}(C,\omega_C) = 0 \quad \Longleftrightarrow \quad p< \mathrm{Cliff}(C)
$$
for every smooth projective curve $C$ over the complex numbers \cite{Green}. Here, $K_{i,j}(C,\omega_C)$ denotes the Koszul cohomology group, see, e.g., \cite{AproduNagelKoszulBook}. Hirschowitz--Ramanan \cite{HirschowitzRamanan} and Voisin \cite{VoisinGreen1,VoisinGreen2} established Green's conjecture for smooth curves lying on K3 surfaces of Picard rank~$1$. Aprodu \cite{AproduSyzygies} subsequently proved that Green's conjecture is implied by the linear growth condition on Brill--Noether loci of pencils,
$$
\dim W^1_d(C) \leq d-\mathrm{gon}(C), \quad \gon(C) \leq d \leq g - \gon(C) + 2.
$$
Aprodu--Farkas then used the linear growth condition to extend Voisin's results to all smooth curves lying on an arbitrary K3 surface \cite{FarkasAproduGreenConj}. Using a similar approach, Lelli-Chiesa established Green's conjecture for smooth curves $C$ on rational surfaces $X$ having an anticanonical pencil, provided $\O_X(C+K_X)$ is big and nef and that, in case $h^0(X,\omega_X^{\vee}) = \chi(X,\omega_X^{\vee})= 2$, the restriction $\omega_X^{\vee}|_C$ does not compute $\Cliff(C)$ for the general curve in $|\O_X(C)|$ \cite{Lelli-ChiesaGreenConj}. Moretti has shown Green's conjecture for smooth curves general in their linear system on abelian surfaces of Picard rank~$1$ with a non-principal polarization \cite{Federico}. Recently, the authors proved that Green's conjecture holds for the general smooth curve in a linear system on a regular surface with nef and effective anticanonical class \cite{GeigerKretschmerBN}.

\noindent In this paper, we extend Lelli-Chiesa's result as well as the applications presented in \cite{GeigerKretschmerBN} and prove Green's conjecture for smooth curves $C$ on arbitrary rational surfaces carrying an effective anticanonical divisor. In fact, keeping with the philosophy of \cite{GeigerKretschmerBN}, it is enough to require that $-K_X|_C$ be effective on $C$. In addition, on birationally ruled surfaces $X$ of irregularity $q \geq 2$ carrying an effective anticanonical divisor, we show that Green's conjecture holds for all smooth curves $C$ on $X$ covering the base curve with degree at least $2$, provided the latter has only finitely many gonal pencils.

\begin{theorem}\label{thm:intro}
    Let $X$ be a smooth projective surface with $-K_X>0$ and $C \subset X$ a smooth curve. Denote by $c$ the Clifford index of the general smooth curve in $|\O_X(C)|$. If $q > 0$, let $f \colon X \to B$ be the ruling and assume $n \coloneqq \deg(f|_C) \geq 2$. Then Green's conjecture holds for $C$ if one of the following conditions holds:
    \begin{enumerate}[(i)]
        \item $q = 0$;
        \item $q \geq 1$ and $n  \gon(B) \neq c +2$;
        \item $q \geq 2$, $n \gon(B)= c+2$ and $B$ has only finitely many gonal pencils.
    \end{enumerate}
\end{theorem}

\noindent Here, $n \gon(B) \neq c+2$ can only happen if $\gon(C) < n \gon(B)$ or $C$ has Clifford dimension at least $2$.

\begin{theorem}[{Corollary~\ref{cor:maincorollary}}]
\label{thm:introRational}
    Let $X$ be a rational surface and $C \subset X$ a smooth curve with $-K_X|_C\geq 0$. Then Green's conjecture holds for $C$. In particular, Green's conjecture holds for all smooth curves on the blow-up of a del Pezzo surface $X$ in $\delta \leq K_X^2$ arbitrary points, possibly infinitely near.
\end{theorem}

\noindent See Corollary~\ref{cor:maincorollary} for a statement similar to Theorem~\ref{thm:introRational} in positive irregularity. As every smooth curve has a nodal plane model and therefore lies on some rational surface, knowing Green's conjecture for smooth curves on arbitrary rational surfaces would confirm the conjecture in general. However, one day prior to the arXiv submission of this paper the preprint \cite{AproduFarkasParkCounterexample} appeared, announcing counterexamples to Green's conjecture that admit a $6$-nodal model given by a cubic section of a quintic del Pezzo surface $Y$. Therefore, counterexamples to Green's conjecture can be realized already on a blowup $X$ of $Y$ in $6 = K_Y^2 + 1$ points! In this sense, Theorem~\ref{thm:introRational} is actually optimal.
\smallskip

\noindent In the case of positive irregularity $q > 0$, Theorem~\ref{thm:intro} overlaps with results of Ballico--Fontanari \cite{BallicoFontanari}, showing Green's conjecture for covers of base curves having only finitely many gonal pencils once the genus of the covering curve $C$ is sufficiently large. However, our results go beyond the hypotheses of Ballico--Fontanari, allowing for genera significantly below their required bound and needing no simplicity hypothesis. On the other hand, Ballico--Fontanari require no embedding of $C$ into an anticanonical surface $X$ ruled over the base curve.

\medskip

\noindent The technique of Lazarsfeld--Mukai bundles, pioneered in \cite{Lazarsfeld}, as well as the work of Donagi--Morrison \cite{DonagiMorrison} in the case of pencils, have proven very useful for studying the Brill--Noether theory of curves on surfaces. To prove Theorem~\ref{thm:intro}, we use a slightly adapted version of Aprodu's linear growth condition and a new parameter count for Lazarsfeld--Mukai bundles. The key difference to previous works is the reduction step in Proposition~\ref{prop:mainreduction}, reducing to the case of curves $C$ disjoint from a fixed effective anticanonical divisor. This reduction allows us to utilize the deformation theory of framed sheaves developed in \cite{BruzzoMarkushevishFramedBundles} to obtain better parameter counts for Lazarsfeld--Mukai bundles. We distinguish the three cases of pencils $(C,A)$ whose Lazarsfeld--Mukai bundles are (i) simple, (ii) non-simple but indecomposable, or (iii) decomposable, i.e., a direct sum of two line bundles. In each case we exploit the geometry of an anticanonical divisor disjoint from $C$ to obtain dimension counts that establish the linear growth condition.

\subsection{Acknowledgements}
S.G. is grateful to Angela Ortega for her encouragement and supervision. S.G. was funded by the Deutsche Forschungsgemeinschaft (DFG, German Research Foundation) — RTG2965 — Project number 512730679 at the Humboldt-University of Berlin. A.K. was supported by the ERC Advanced Grant SYZYGY (No. 834172, PI: Gavril Farkas).

\subsection{Statement on the use of AI}
The basis for this project stems from our previous work \cite{GeigerKretschmerBN} and our intuition that its assumption of irregularity zero may be circumvented for the application to Green's conjecture. The output of an LLM, prompted with this intuition as its input, provided a working basis for a rigorous proof. The idea of separating the curve from the anticanonical divisor by a sequence of blow-ups was part of an earlier conversation with the LLM in the rational case. In this sense, a significant part of the paper relies on the AI contribution. The authors take full responsibility for the veracity of the results and proofs presented here.

\section{Proofs}

\subsection{Linear growth and reduction to a disjoint anticanonical divisor}
 Let $V$ be a complex vector space and $M$ a graded module over the symmetric algebra $S^*V$. The \emph{Koszul cohomology group $K_{p,q}(M,V)$} is the cohomology at the middle term of the complex
 $$
 \bigwedge^{p+1}V\otimes M_{q-1}\to \bigwedge^pV\otimes M_{q}\to \bigwedge^{p-1}V\otimes M_{q+1}, 
 $$
 where the maps are induced by the dual of the wedge product map. If $X$ is a projective variety and $L$ a line bundle on $X$, then the Koszul cohomology of $L$, denoted by $K_{p,q}(X,L)$, is the Koszul cohomology of $R(L)\coloneqq \bigoplus_iH^0(X,L^i)$ considered as a graded module over $S^*\!H^0(X,L)$. For a smooth curve $C$, Green's conjecture predicts the vanishing
 \begin{equation*}
    K_{p,2}(C,\omega_C) = 0 \quad \Longleftrightarrow \quad p < \Cliff(C).
\end{equation*}
The conjectural part reduces to $K_{c-1,2}(C,\omega_C) = 0$ or, equivalently, $K_{g-c-1,1}(C,\omega_C) = 0$. As stated in the introduction, Aprodu \cite{AproduSyzygies} proved that Green's conjecture for a curve $C$ of genus $g$ and gonality $k$ is implied by the growth condition
\begin{equation*}
    \dim W^1_{d}(C) \leq d-k, \quad \text{for all} \ \ k\leq d \leq g-k+2.
\end{equation*}
Aprodu and Farkas provide a slightly weaker linear growth condition in \cite[Theorem~4.5]{AproduFarkasKoszulCohomology}. As explained in \cite[Section 5.4]{GeigerKretschmerBN}, except for the case of odd genus and maximal gonality which has been settled by Hirschowitz--Ramanan \cite{HirschowitzRamanan} and Voisin \cite{VoisinGreen2}, Green's conjecture is implied by
\begin{equation}\label{eqn:lineargrowthcondition}
    \dim \widetilde{W}^{1,\mathrm{bpf}}_{d}(C) \leq d - c - 2, \quad \text{for all} \ \ c+2 \leq d \leq g-c,
\end{equation}
where $c = \Cliff(C)$ is the Clifford index of $C$ and $\widetilde{W}^{1,\mathrm{bpf}}_d(C)$ denotes the Brill--Noether locus of basepoint-free and complete pencils on $C$ of degree $d$. We will refer to \eqref{eqn:lineargrowthcondition} as the \emph{linear growth condition} even though it deviates slightly from Aprodu's classical formulation. The advantage is that \eqref{eqn:lineargrowthcondition} allows us to treat all Clifford dimensions equally.

\noindent On surfaces with effective anticanonical divisor, we have the following crucial reduction step:

\begin{proposition}[Reduction to a disjoint anticanonical divisor]\label{prop:mainreduction}
The following are equivalent:
\begin{enumerate}
    \item For all smooth projective surfaces $X$ with $-K_X \in \Pic(X)$ effective and nonzero, Green's conjecture holds for all smooth curves $C \subset X$ that are not contained in the base locus of $|-K_X|$.
    \item For all smooth projective surfaces $X$ with $-K_X \in \Pic(X)$ effective and nonzero, every fixed anticanonical divisor $E \in |-K_X|$ and any $\mathcal{L} \in \Pic(X)$, Green's conjecture holds for the general smooth curve $C \in |\mathcal{L}|$ that is disjoint from $E$.
\end{enumerate}
\end{proposition}

\begin{proof}
    The implication $(1) \Rightarrow (2)$ is clear.
    For the reverse implication, let $C \subset X$ be a smooth curve not contained in the base locus of $|-K_X|$. Then there is an anticanonical divisor $E \in |-K_X|$ not containing $C$. If $C$ is not already disjoint from $E$, denote by $\pi \colon X'\to X$ the consecutive blowup of the reduced intersection points of $C$ and $E$. Denote by $F_i$ the exceptional divisors and set $C' \coloneqq \pi^*C - \sum_i F_i$ and $E' \coloneqq \pi^*E - \sum_i F_i \in |-K_{X'}|$; in particular, $E'$ is still an anticanonical divisor. We have $C'\cap E'=\emptyset$ and, as $C$ is smooth, $C \simeq C'$. Therefore, Green's conjecture holds for $C$ if and only if it holds for $C'$. We now show that Green's conjecture for $C'$ does not depend on the choice of curve in its linear system on $X'$: As $C'$ and $E'$ are disjoint, we have 
    \begin{equation*}
        \O_{X'}(C')|_{E'}\simeq \O_{E'} \quad \text{ and } \quad \O_{X'}(C')|_{C'} \simeq \omega_{C'} .
    \end{equation*}
    Hence, $C'^2 >0$ and $C'$ is big and nef. Consider the evaluation sequence
    $$
    0\to \O_{X'}(K_{X'}+mC') \to \O_{X'}(mC')\to \O_{X'}(mC')|_{E'}\simeq \O_{E'}\to 0.
    $$
    By Kawamata--Viehweg vanishing, $H^i(X',\O_{X'}(K_{X'}+mC'))= 0$ for $i > 0$, so we get the isomorphisms 
    $$
    H^1(X',\O_{X'}(mC')) \xrightarrow{\sim} H^1(E',\O_{E'})
    $$
    for all $m \geq 1$. Let $s_{C'}$ be the section defining $C'$ and consider the restriction to $E'$ of multiplication by $s_{C'}$. We obtain the commutative diagram
    \[\begin{tikzcd}
	   {\mathcal{O}_{X'}((m-1)C')} & {\mathcal{O}_{X'}(mC')} \\
	   {\mathcal{O}_{E'}} & {\mathcal{O}_{E'}}
	   \arrow["{s_{C'}}", hook, from=1-1, to=1-2]
	   \arrow[two heads, from=1-1, to=2-1]
	   \arrow[two heads, from=1-2, to=2-2]
	   \arrow["\simeq"', from=2-1, to=2-2].
    \end{tikzcd}\]
     The lower horizontal map is an isomorphism as $s_{C'}$ is everywhere nonzero on $E'$. Moreover, both vertical maps become isomorphisms on $H^1$ for $m\geq 2$, so $s_{C'}$ induces isomorphisms
     \begin{equation}
     \label{eqn:multiplicationiso}
        H^1(X',\O_{X'}((m-1)C'))\xrightarrow{\sim} H^1(X',\O_{X'}(mC')).
     \end{equation}
     For $m=1$, the left vertical map in the above diagram is still injective on $H^1$: This is clear when $h^1(X',\O_X') = 0$. If $q \coloneqq h^1(X',\O_{X'}) \geq 1$, then a minimal model of $X'$ is a geometrically ruled surface over a smooth base curve $B$ of genus $q$, inducing a ruling map $f \colon X' \to B$. Let $\Gamma$ be a component of $E'$ dominating $B$. The finite flat map $f|_{\Gamma} \colon \Gamma \to B$ induces an injection $f|_{\Gamma}^\ast \colon H^1(B,\O_B) \hookrightarrow H^1(\Gamma,\O_{\Gamma})$ because it has a left inverse given by $\frac{1}{\deg(f|_\Gamma)} H^1(\mathrm{Tr})$, where $\mathrm{Tr} \colon  (f|_\Gamma)_\ast \O_\Gamma \to \O_B$ is the trace map. Under the identification $H^1(B,\O_B)\simeq H^1(X',\O_X')$ given by pullback along the ruling, we find that the composition 
     $$
     H^1(X',\O_X')\to H^1(E',\O_{E'})\to H^1(\Gamma,\O_{\Gamma})
     $$
     given by restriction, is exactly $f|_{\Gamma}^\ast$, thus $H^1(X',\O_X') \to H^1(E',\O_E')$ is injective. It follows that \eqref{eqn:multiplicationiso} is still injective for $m = 1$. The evaluation sequence for $s_{C'}$,
    $$
    0\to \O_{X'}((m-1)C')\to \O_{X'}(mC')\to \O_{C'}(mK_{C'})\to 0,
    $$
     then gives surjections
     $$
     H^0(X',\O_{X'}(mC'))\twoheadrightarrow H^0(C',\O_{C'}(mK_{C'}))
     $$
      for all $m\geq 1$ and consequently an isomorphism of graded rings $R(X',\O_{X'}(C'))/(s_{C'})\simeq R(C',\O_C'(K_{C'}))$. As $s_{C'}$ is a (degree~$1$) homogeneous regular element in the graded ring $R(X',\O_{X'}(C'))$, the quotient by $(s_{C'})$ preserves minimal free resolutions and we obtain
      \begin{equation}
          \label{eqn:KoszulEqualityReduction}
          K_{p,q}(C',K_{C'})\simeq K_{p,q}(X',\O_{X'}(C')).
      \end{equation}
      It follows that Green's conjecture for smooth curves disjoint from an anticanonical divisor does not depend on the choice of the curve in its linear system.
\end{proof}

\begin{remark}
    For a smooth curve $C$ of genus $g \geq 2$ on a birationally ruled surface $f \colon X \to B$ with $n \coloneqq \deg(f|_C) \geq 2$, the curve cannot be contained in the base locus of $-K_X$. Indeed, if $-K_X = mC+D$ for some $m>0$ and $D$ effective, let $F$ be the class of a fiber of $f$. We get $2 = -K_X.F = mF.C + F.D \geq m n$ since $F$ is nef. Hence, $m = 1$, $n = 2$ and adjunction yields a contradiction to $g \geq 2$.
\end{remark}

\subsection{Generic Green's conjecture for curves disjoint from an anticanonical divisor} Let $X$ be a smooth projective surface with $-K_X>0$ and $C\subset X$ a smooth curve. If the irregularity $q \coloneqq h^1(X,\O_X)$ is zero, then $X$ is rational. If $q > 0$, a minimal model of $X$ is geometrically ruled, so there is a map $f \colon X \to B$ onto a smooth projective curve $B$ of genus $q$. This morphism is unique up to automorphisms of $B$; it is the Albanese morphism of $X$ with respect to some basepoint. Slightly abusing terminology, we call $f$ \emph{the ruling of $X$}. We fix the following notation:
\begin{itemize}
    \item If $q>0$, then $n\coloneqq \mathrm{deg}(f|_C)$ denotes the degree of the finite map $f|_C\colon C\to B$;
    \item $\mathcal{L} \coloneqq \O_X(C)$ and $|\mathcal{L}|_{\mathrm{sm}}$ is the open locus of smooth, irreducible curves in $|\mathcal{L}|$;
    \item $c$ denotes the Clifford index of the general curve in $|\mathcal{L}|_{\mathrm{sm}}$.
\end{itemize}
By Proposition~\ref{prop:mainreduction}, it suffices to establish generic Green's conjecture for curves disjoint from a fixed effective anticanonical divisor.

\begin{definition}
    Let $X$ be a smooth projective surface and $i\colon C \hookrightarrow X$ a smooth curve. Let $A$ be a basepoint-free line bundle on $C$. We define the \emph{Lazarsfeld--Mukai bundle $F_{C,A}$ associated to $(C,A)$} by the sequence
    $$
    0 \to F_{C,A}^{\vee} \to H^0(C,A) \otimes \O_X \to i_\ast A \to 0.
    $$
    Its dual sequence is given by
    $$
    0 \to H^0(C,A)^{\vee} \otimes \O_X \to F_{C,A} \to i_\ast(A^{\vee} \otimes \O_C(C)) \to 0
    $$
    and $F_{C,A}$ is a vector bundle of rank $\rk(F_{C,A}) = h^0(C,A)$.
\end{definition}

\noindent Fix an anticanonical divisor $E\in |-K_X|$ disjoint from $C$. We want to establish the linear growth condition (\ref{eqn:lineargrowthcondition}) for $\widetilde{W}^{1,\mathrm{bpf}}_{d}(C)$. Let $\mathcal{W} \subseteq \mathcal{W}^1_d(|\mathcal{L}|)$ therefore be an irreducible component dominating $|\mathcal{L}|$ and consisting generically of basepoint-free and complete pencils on a smooth curve, i.e., for the general $(C,A) \in \mathcal{W}$ we have that $C$ is smooth and $A$ is basepoint-free with $h^0(A) = 2$. Equivalently, letting $|\mathcal{L}|_{\mathrm{sm}}$ be the open locus of smooth, irreducible curves in $|\mathcal{L}|$ and $\widetilde{\mathcal{W}}^{1,\mathrm{bpf}}_d(|\mathcal{L}|_{\mathrm{sm}})$ the relative Brill--Noether locus of basepoint-free and complete pencils, we may take $\mathcal{W}$ to be a dominating irreducible component of $\widetilde{\mathcal{W}}^{1,\mathrm{bpf}}_d(|\mathcal{L}|_{\mathrm{sm}})$ directly.
Let $F_{C,A}$ be the Lazarsfeld--Mukai bundle associated to $(C,A) \in \mathcal{W}$. Three cases can occur for the general $(C,A) \in \mathcal{W}$:
\begin{enumerate}[(i)]
    \item $F_{C,A}$ is simple;
    \item $F_{C,A}$ is non-simple but indecomposable;
    \item $F_{C,A}$ is decomposable, i.e., a direct sum of two line bundles.
\end{enumerate}
For $q= 0$, one may use \cite[Proposition~3.13]{GeigerKretschmerBN} in order to immediately establish linear growth in case (i) but this is no longer available for positive irregularity. We will instead use the deformation theory of framed sheaves \cite{BruzzoMarkushevishFramedBundles} which is applicable because, via Proposition~\ref{prop:mainreduction}, we have the framing $F_{C,A}|_E \simeq \O_E^{\oplus 2}$.
We start with the following observation.

\begin{lemma}
\label{lem:constancyofLMsections}
    Let $X$ be a surface with $-K_X>0$ and $E\in |-K_X|$ a fixed anticanonical divisor. Let $C\subset X$ be a smooth curve of genus $g$ disjoint from $E$ and $A$ a complete and basepoint-free pencil of degree $d$ on $C$. Then $h^0(X,F_{C,A}) = g-d+3$ and $h^1(X,F_{C,A}) = 2q+2$.
\end{lemma}

\begin{proof}
    Consider the defining Lazarsfeld--Mukai sequence
    $$
    0 \to F_{C,A}^{\vee}\to H^0(C,A)\otimes \O_X\to A\to 0.
    $$
    As $C\cap E = \emptyset$, we again have $\O_X(C)|_E \simeq \O_E$ and $\O_X(C)|_C \simeq \omega_C$. Let $s_E\in H^0(X,\omega_X^{\vee})$ be the section defining $E$. Dualizing the above sequence, we obtain the commutative diagram

    \[\begin{tikzcd}
	   0 & {\omega_X^{\oplus 2}} &     {F_{C,A}\otimes\omega_X} & {A^{\vee}\otimes\omega_C} & 0 \\
	   0  & {\mathcal{O}_X^{\oplus 2}} & {F_{C,A}} & {A^{\vee}\otimes\omega_C} & 0
	   \arrow[from=1-1, to=1-2]
	   \arrow[from=1-2, to=1-3]
	   \arrow["{s_E}", from=1-2, to=2-2]
	   \arrow[from=1-3, to=1-4]
	   \arrow["{s_E}", from=1-3, to=2-3]
	   \arrow[from=1-4, to=1-5]
	   \arrow["\sim", from=1-4, to=2-4]
	   \arrow[from=2-1, to=2-2]
	   \arrow[from=2-2, to=2-3]
	   \arrow[from=2-3, to=2-4]
	   \arrow[from=2-4, to=2-5].
    \end{tikzcd}\]
    The right vertical map is an isomorphism as $s_E$ is nonvanishing on $C$ by definition. Taking cohomology of the above diagram, we can see that, if $H^1(s_E)\colon H^1(X,\omega_X)\to H^1(X,\O_X)$ is the zero map, then so is the connecting morphism $\delta\colon H^0(C,A^{\vee}\otimes \omega_C)\to H^1(X,\O_X)^{\oplus 2}$ of the bottom exact sequence. Then global sections are additive in that sequence and thus $h^0(X,F_{C,A}) = g-d+3$ and $h^1(X,F_{C,A}) = 2q+2$. Now, $H^1(s_E) = 0$ is clear if $h^1(X,\O_X) = 0$, and if $h^1(X,\O_X) >0$, the anticanonical evaluation sequence
    $$
    0\to \omega_X\xrightarrow{s_E}\O_X \to \O_E \to 0
    $$
    shows that we only need the injectivity of $H^1(X,\O_X) \to H^1(E,\O_E)$. This in turn we have shown in the proof of Proposition~\ref{prop:mainreduction}.
\end{proof}

\begin{proposition}
\label{prop:genericgreen1}
   Let $\mathcal{W} \subseteq \widetilde{\mathcal{W}}^{1,\mathrm{bpf}}_d(|\mathcal{L}|_{\mathrm{sm}})$ be a dominating irreducible component as above such that, for the general $(C,A)\in \mathcal{W}$, the Lazarsfeld--Mukai bundle $F_{C,A}$ is simple. Then the general fiber of $\mathcal{W}$ over $|\mathcal{L}|$ has the expected dimension $\rho(g,1,d)$ and in particular satisfies the linear growth condition \eqref{eqn:lineargrowthcondition}.
\end{proposition}

\begin{proof}
    Let $(C,A)\in \mathcal{W}$ be a basepoint-free complete pencil such that $F= F_{C,A}$ is simple. Let $i \colon C \hookrightarrow X$ be the inclusion and consider the Lazarsfeld--Mukai exact sequence
    $$
    0\to F^{\vee}\to H^0(C,A)\otimes \O_X\to i_*A\to 0.
    $$
    As $i_*A$ is supported on $C$ and $C$ is disjoint from $E$, dualizing and restricting to $E$ yields an isomorphism $F|_E\simeq \O_E^2$. The set of such Lazarsfeld--Mukai bundles is thus dominated by the parameter space of simple bundles whose restriction to $E$ is trivial. Fix an isomorphism $F|_E\simeq \O_E^2$ and consider the pair $(F,E)$. The deformation theory of framed bundles \cite[Section 4]{BruzzoMarkushevishFramedBundles} applies verbatim to simple sheaves and the tangent space at $(F,E)$ in the moduli space of framed bundles is given by $\Ext^1(F,F(K_X))$. An obstruction space is given by the trace-free part of $\Ext^2(F,F(K_X))$, denoted by $\Ext^2(F,F(K_X))_0 \subseteq \Ext^2(F,F(K_X))$. We have $\Ext^2(F,F(K_X)) \simeq \mathbb{C}$ and the trace map $\mathrm{tr}\colon \Ext^2(F,F(K_X))\to H^2(X,\omega_X)$ is nonzero, thus the moduli space of framed simple bundles is unobstructed.  As $F$ is simple and $-K_X > 0$, we have $\Hom(F,F(K_X)) = 0$. A Hirzebruch--Riemann--Roch calculation then shows
    $$
    \dim \Ext^1(F,F(K_X)) = -2g-1+4c_2(F)+4h^1(X,\O_X).
    $$
    If $X$ is regular, then $h^0(E,\O_E) = 1$ for all anticanonical divisors $E$. If $q > 0$, denote by $f \colon X \to B$ the ruling. Pushing forward the short exact sequence
    \begin{equation*}
        0 \rightarrow \omega_X \rightarrow \O_X \rightarrow \O_E \rightarrow 0
    \end{equation*}
    along $f$ yields
    \begin{equation*}
        0 \rightarrow \O_B \rightarrow f_\ast \O_E \rightarrow R^1 f_\ast \omega_X \rightarrow 0
    \end{equation*}
    since, by Grauert's theorem, $f_\ast \omega_X = 0$, $f_\ast \O_X = \O_B$ and $R^1 f_\ast \O_X = 0$. Moreover, $\omega_{X/B} \simeq \omega_X \otimes f^\ast(\omega_B)^\vee$ by the smoothness of $X$ and $B$ and using relative adjunction for the factorization of $f$ as $X \rightarrow X \times B \rightarrow B$ given by the graph of $f$, embedding $X$ into $X \times B$ as a Cartier divisor, followed by projection to $B$. Relative duality then gives $R^1 f_\ast \omega_X = R^1 f_\ast (\omega_{X/B} \otimes f^\ast(\omega_B)) = R^1 f_\ast(\omega_{X/B}) \otimes \omega_B \simeq \omega_B$ since $R^1 f_\ast(\omega_{X/B}) \simeq \mathcal{H}om(f_\ast(\O_X), \O_B) = \mathcal{H}om(\O_B, \O_B) = \O_B$.
    We therefore obtain the short exact sequence
    $$
    0 \to \O_B \to f_\ast \O_E \to \omega_B \to 0.
    $$
    The sequence shows that $f_\ast \O_E$ is locally free of rank $2$, so the trace map $\frac{1}{2}\mathrm{Tr}\colon f_*\O_E\to \O_B$ gives a left inverse to the inclusion $\O_B \rightarrow f_\ast \O_E$ and therefore splits the sequence. Hence, $h^0(E,\O_E) = 1+ g(B) = 1 + h^1(X,\O_X)$. The above formula becomes
    $$
    \dim \Ext^1(F,F(K_X)) = -2g-5+4c_2(F)+4h^0(E,\O_E).
    $$
    Passing from framed bundles to isomorphism classes of Lazarsfeld--Mukai bundles, every fiber has dimension at least $\dim \mathrm{Aut}(\O_E^2)-1 = 4h^0(E,\O_E)-1$ because the automorphisms of $\O_E^2$ act faithfully on the framing isomorphisms.
    As $\mathcal{W}$ is in turn dominated by the relative Grassmannian of $2$-dimensional subspaces of the global sections over the set of relevant Lazarsfeld--Mukai bundles, Lemma \ref{lem:constancyofLMsections} and $d=c_2(F)$ gives
    \begin{align*}
        \dim \mathcal{W}&\leq \dim \Ext^1(F,F(K_X)) - (\dim \Aut(\O_E^2)-1) + \dim \Gr(2,H^0(F)) \\
        &= -2g-5+4d+4h^0(E,\O_E) - (4h^0(E,\O_E)-1) + 2(g-d+1) \\
        &= 2d-2 = \dim |\mathcal{L}| + 2d-2 -g = \dim |\mathcal{L}| + \rho(g,1,d).
    \end{align*}
    Finally, as we only need to check the linear growth condition \eqref{eqn:lineargrowthcondition} for $d \leq g-c$, this shows
    $$
    \dim \mathcal{W}\leq \dim |\mathcal{L}| +d-c-2,
    $$
    as desired.
\end{proof}

\begin{proposition}[{\cite[Lemma~4.4]{DonagiMorrison}, \cite[Lemma~2.12]{GeigerKretschmerBN}}]\label{prop:DMextension}
    Let $X$ be a smooth projective surface and $F$ a generically globally generated rank~$2$ vector bundle on $X$ which is neither simple nor decomposable. Then $F$ sits in a \emph{Donagi--Morrison extension}
    \begin{equation*}
        0 \rightarrow \O_X(H+H') \rightarrow F \rightarrow \mathcal{I}_Z(H) \rightarrow 0,
    \end{equation*}
    where $H, H' \in \Pic(X)$ are effective divisor classes with $H$ nonzero, $2H+H' = c_1(F)$, and $\mathcal{I}_Z$ is the ideal sheaf of a zero-dimensional lci subscheme $Z \subset X$, possibly empty.
\end{proposition}

\begin{remark}
    We will make repeated use of the fact that, when in Proposition~\ref{prop:DMextension} $F = F_{C,A}$ is a Lazarsfeld--Mukai bundle, then $H$ has no fixed curves. This follows from the integrality of $C$.
\end{remark}

\begin{lemma}\label{lem:DimensionofTinequality}
    Let $\mathcal{W} \subseteq \widetilde{\mathcal{W}}^{1,\mathrm{bpf}}_d(|\mathcal{L}|_{\mathrm{sm}})$ be a dominating irreducible component such that, for the general $(C,A) \in \mathcal{W}$, the Lazarsfeld--Mukai bundle $F_{C,A}$ is non-simple but indecomposable. For fixed $l$, let $\mathcal{T} \subseteq \mathrm{Hilb}^l(X) \times \Pic(X)$ be an irreducible, locally closed stratum of pairs $(Z,H)$ occurring in the Donagi--Morrison extension of some $F_{C,A}$, $(C,A) \in \mathcal{W}$, such that the numerical class of $H$ as well as $h^1(E,\O_E(H))$ and $h^1(E,\mathcal{I}_Z(H)|_E)$ are constant along $\mathcal{T}$. Define $\mathcal{P} \coloneqq \mathrm{im}(\mathcal{T}\to \Pic(X))$ and $p \coloneqq \dim \mathcal{P}$. Then
    $$
    \dim \mathcal{T} \leq p + 2l + 2(h^1(E,\O_E(H)) - h^1(E,\mathcal{I}_Z(H)|_E)).
    $$
\end{lemma}

\begin{proof}
      Clearly, we only need to bound the dimension of the fiber $\mathcal{T}_{H}$ of $\mathcal{T} \to \Pic(X)$ over $\O_X(H)$. Restricting the Donagi--Morrison exact sequence gives
      $$
        0\to \O_X(-H)|_E \to \O_E^2 \to \mathcal{I}_Z(H)|_E \to 0
      $$
      because $\O_X(C)|_E = \O_E$ and $2H+H' = [C] \in \Pic(X)$. Left-exactness follows because $\mathcal{T}or_1^{\O_X}(\mathcal{I}_Z(H), \O_E) = \ker(\mathcal{I}_Z(H-E) \rightarrow \mathcal{I}_Z(H)) = 0$. Consider
      \begin{align*}
        \mathcal{U}_H &\coloneqq \{(Z,\alpha)\colon Z\in \mathcal{T}_{H}, \alpha \colon \O_E^2 \twoheadrightarrow \mathcal{I}_Z(H)|_E \}, \\ \mathcal{K}_H &\coloneqq \{K \subseteq \O_E^2 \colon \O_E^2/K \simeq \mathcal{I}_Z(H)|_E \text{ for some } Z \in \mathcal{T}_H\}
     \end{align*}
     and the map $\pi \colon \mathcal{U}_H \to \mathcal{K}_H, \ (Z,\alpha) \mapsto \ker(\alpha)$. Abstractly, $\ker(\alpha) \simeq \O_X(-H)|_E$ by a variant of Schanuel's lemma: Forming the fiber product $P \subseteq \O_E^2 \oplus \O_E^2$ of the two surjections $\O_E^2 \twoheadrightarrow \mathcal{I}_Z(H)|_E$ given by $\alpha$ and the surjection in the restricted Donagi--Morrison sequence shows, first, that $P$ is locally free of rank~$3$ since the kernel of the Donagi--Morrison sequence is a line bundle, and second, that $\ker(\alpha)$ must then also be a line bundle on $E$ because it is the kernel of a surjection $P \twoheadrightarrow \O_E^2$. Taking determinants then gives $\ker(\alpha) \simeq \O_E(-H)$. We get
   \begin{align*}
   \dim \mathcal{K}_H &\leq \dim \Hom(\O_E(-H), \O_E^2) - \dim \mathrm{Aut}(\O_E(-H)) \\
   &= 2h^0(E,\O_E(H))-(h^1(X,\O_X)+1) \\
   &= 2(h^1(E,\O_E(H))-K_X.H)-h^1(X,\O_X)-1.
   \end{align*}
    Let $(K \subseteq \O_E^2) \in \mathcal{K}_H$ be fixed and set $Q \coloneqq \O_E^2/K$. The fiber $\pi^{-1}(K)$ then parametrizes pairs of subschemes $Z$ and isomorphisms $\mathcal{I}_Z(H)|_E \xrightarrow{\sim} Q$. Up to quotienting by scalar automorphisms, this is again the data of a framing of $\mathcal{I}_Z(H)$ along $E$. Since $\mathcal{I}_Z(H)$ is torsion-free, the arguments in \cite[Section 4]{BruzzoMarkushevishFramedBundles} apply verbatim and we retain the dimension bound for $\pi^{-1}(K)$ given by
    $$
    \dim \pi^{-1}(K)\leq \ext^1(\mathcal{I}_Z(H),\mathcal{I}_Z(H-E))+1 = h^1(X,\O_X)+2l+1.
    $$
    The last equation uses Riemann--Roch together with $\hom(\mathcal{I}_Z(H),\mathcal{I}_Z(H-E)) = 0$ and
    \begin{equation*}
        \ext^2(\mathcal{I}_Z(H),\mathcal{I}_Z(H-E)) = \hom(\mathcal{I}_Z(H),\mathcal{I}_Z(H)) = 1.
    \end{equation*}
    Every fiber of the map $\mathcal{U}_H \rightarrow \mathcal{T}_H$ has dimension $\hom(\O_E^2, \mathcal{I}_Z(H)|_E)$ because the general homomorphism is surjective once a surjection exists. Therefore, 
    \begin{equation*}
        \dim \mathcal{T}_{H} = \dim \mathcal{U}_H- \dim \Hom(\O_E^2, \mathcal{I}_Z(H)|_E) = \dim \mathcal{U}_H - 2h^0(E,\mathcal{I}_Z(H)|_E).
    \end{equation*}
    Using $h^0(E,\mathcal{I}_Z(H)|_E) = h^1(E,\mathcal{I}_Z(H)|_E)-K_X.H$, we get
    \begin{align*}
    \dim \mathcal{T}_{H} &\leq \dim \mathcal{K}_H+\mathrm{sup}_K\dim \pi^{-1}(K)-2h^0(E,\mathcal{I}_Z(H)|_E) \\
    &\leq 2l+2(h^1(E,\O_E(H))-h^1(E,\mathcal{I}_Z(H)|_E)).\qedhere
    \end{align*}
\end{proof}

\begin{lemma}\label{lem:PicardBound}
    With the same notation as in Lemma~\ref{lem:DimensionofTinequality}, assume that $-K_X.H>0$ for $H \in \mathcal{P}$. Then
    $$
    p \leq q-h^1(E,\O_E(H)).
    $$
\end{lemma}

\begin{proof}
    Let $s\in H^0(E,\O_E(H))$ be a nonzero regular section and consider
    $$
    \alpha_s\colon H^1(X,\O_X)\to H^1(E,\O_E)\xrightarrow{s}H^1(E,\O_E(H)).
    $$
    In the proof of Proposition \ref{prop:mainreduction}, we have shown that the first map is injective. Using Serre duality and $\omega_E \simeq \O_E$, we find that the dual map $\alpha_s^{\vee}$ is given by
    $$
    \alpha_s^{\vee} \colon H^0(E,\O_E(-H)) \xrightarrow{\cdot s} H^0(E,\O_E)/(\C\cdot 1).
    $$
    To see this, consider the anticanonical divisor sequence
    $$
    0\to \O_X(K_X)\to \O_X\to \O_E\to 0.
    $$
    As restriction $H^1(X,\O_X)\to H^1(E,\O_E)$ is injective, the boundary map $\delta\colon H^0(E,\O_E)\to H^1(X,\O_X(K_X))$ is surjective with kernel given by the image of $1$ under $H^0(X,\O_X)\to H^0(E,\O_E)$. Serre duality thus identifies
    $$
    H^1(X,\O_X)^{\vee}\simeq H^1(X,\O_X(K_X)) \simeq H^0(E,\O_E)/(\mathbb{C}\cdot 1).
    $$
     Notice that $s$ is regular and $E.H = -K_X.H>0$ implies that $s$ does not vanish on all of $E$, thus $\alpha_s^{\vee}$ is injective and $\alpha_s$ is surjective. Therefore, $\dim \ker(\alpha_s) = q-h^1(E,\O_E(H))$. By the assumptions on $\mathcal{T}$, the number $h^0(E,\O_E(H))$ is constant for $H \in \mathcal{P}$. Therefore, the section $s$ lifts along all first order deformations of a fixed $H$ \emph{along} $\mathcal{P}$, essentially by Grauert's theorem. In particular, if $\xi \in H^1(X,\O_X)$ is a tangent direction along which $s$ lifts, then $\alpha_s(\xi) = 0$ by \cite[Proposition~3.3.4]{Sernesi}, so we obtain $T_H \mathcal{P} \subseteq \ker(\alpha_s)$. In particular, $p \leq q-h^1(E,\O_E(H))$. (Technically, one should pass from the constructible set $\mathcal{P}$ to its locally closed constituents and then apply this argument.)
\end{proof}

\begin{lemma}\label{lem:pullbackreduction}
    With the same notation as in Lemma \ref{lem:DimensionofTinequality}, assume $q \geq 2$ and let $f \colon X \to B$ be the ruling. Let $F$ be the numerical class of a fiber, and $n = \deg(f|_C) = C.F$. Assume $H.F > 0$ for $H \in \mathcal{P}$. Then
    $$
    n \gon(B) - 3 \leq d-l + K_X.H.
    $$
\end{lemma}

\begin{proof}
    The general fiber of $f\colon X\to B$ is isomorphic to $\P^1$. Set $n = C.F$, then $\mathcal{L}|_F \simeq \O_{\P^1}(n)$. The complete linear system associated to $\O_{\P^1}(n)$ is linearly semistable and $\mathcal{L} = \O_X(C)$ is nef on $X$, so the assumptions of \cite[Proposition 3.1]{BarjaStoppinoSlopeinequality} are satisfied. We obtain the inequality
    $$
    C^2\geq \frac{2n}{n+1}\deg(f_\ast \O_X(C)).
    $$
    As effective line bundles on trees of smooth rational curves have vanishing $H^1$, Grauert's theorem gives $R^1 f_\ast \O_X(C) = 0$, and $R^2 f_\ast \O_X(C) = 0$ for dimension reasons. Therefore, $f_\ast \O_X(C)$ is locally free of rank $n+1$, and Leray's spectral sequence stabilizes on the second page, implying $\chi(B,f_*\O_X(C)) = \chi(X,\O_X(C))$. Both sides may be computed by Riemann--Roch, giving
    \begin{equation*}
        C^2 = 2\deg(f_\ast \O_X(C)) - 2n(q-1) + C.K_X.
    \end{equation*}
    Combining this equation with the above inequality gives
    \begin{equation*}
        \frac{2}{n+1} \deg(f_\ast \O_X(C)) + C.K_X \geq 2n(q-1).
    \end{equation*}
    Multiplying this inequality by $n$, we obtain 
    $$
    C^2 + n C.K_X \geq \frac{2n}{n+1} \deg(f_\ast \O_X(C)) + n C.K_X \geq 2n^2 (q-1).
    $$
    The same inequality holds for $H$ in place of $C$ because $H.F > 0$ and $R^i f_\ast \O_X(H) = 0$ for $i > 0$ since the restriction of $H$ to every fiber of $f$ is effective as $H$ has no fixed curves by construction. In the end, since $C.K_X = 0$, we obtain
    \begin{equation}
    \label{eqn:lemmaproofeqn1}
    C^2 \geq n^2(2q-2) \quad \text{ and } \quad H^2+\delta\cdot H.K_X\geq \delta^2(2q-2),
    \end{equation}
    where $\delta \coloneqq H.F > 0$. The class $H-\frac{\delta}{n}C$ is orthogonal to the nef class $F$, so the Hodge index theorem implies
    $$
    (H-\frac{\delta}{n}C)^2\leq 0 \quad \text{ or, equivalently, } \quad H^2+\frac{\delta^2}{n^2}C^2\leq \frac{2\delta}{n}H.C.
    $$
    Since $F$ is nef and $[C] = 2H + H'$, we have $n \geq 2 \delta$. Consequently, using $l = d - H.(H+H')$, we get
    \begin{align}
    \label{eqn:lemmaproofeqn2}
    \begin{split}
        d-l+K_X.H &= H.C-H^2+K_X.H\\
        &\geq \left(\frac{n}{2 \delta} - 1\right) H^2 + \frac{\delta}{2n} C^2 + H.K_X\\
        &\geq \left(\frac{n}{2 \delta} - 1\right) \left( \delta^2(2q-2) - \delta H.K_X\right) + \frac{\delta}{2n} n^2(2q-2) + H.K_X\\
        &= (2q-2)\delta(n-\delta) + \frac{n - 2\delta - 2}{2} H.(-K_X).
    \end{split}
    \end{align}
    Since $\delta \geq 1$ and $n \geq 2\delta$, we always have $\delta(n-\delta) \geq n-1$. As $H$ has no fixed curve, $H$ is nef. In case $n - 2 \delta \geq 2$, since $-K_X$ is effective, we get
    $$
    d-l+K_X.H\geq (2q-2)(\delta n-\delta^2)\geq (2q-2)(n-1)\geq nq-2\geq n \mathrm{gon}(B)-3,
    $$
    as desired.
    Let now $n-2\delta\in\{0,1\}$. By \eqref{eqn:lemmaproofeqn1}, nefness of $H$ and effectivity of $H'$ and $-K_X$, we have
    $$
    d-l+K_X.H = H.(H+H') + H.K_X \geq H^2 + \delta H.K_X \geq \delta^2 (2q-2).
    $$
    As $n-2\delta\in \{0,1\}$ we have $n\leq 2\delta+1$. If $\delta \geq 2$, we obtain $\delta^2(2q-2) \geq 2 \frac{n-1}{2} (2q-2) \geq nq-3$, and the same is true for $(\delta,n) = (1,2)$. The only remaining case is $(\delta,n) = (1,3)$. Here, we have $H'.F = n - 2 \delta = 1$ and hence $(K_X+H').F = -1$. This implies $f_\ast \O_X(K_X+H') = 0$ because any local section has to vanish at the non-empty open subset over which $f$ is a $\P^1$-bundle and hence has to vanish identically. Moreover, the sheaf $R^1 f_\ast \O_X(K_X+H')$ is torsion as its general fiber is $H^1(\P^1,\O_{\P^1}(-1)) = 0$. As $R^2 f_\ast \O_X(K_X+H') = 0$ for dimension reasons, the Leray spectral sequence again stabilizes at the second page, implying
    $$
    \chi(X,\O_X(K_X+H')) = -\mathrm{length}(R^1 f_\ast \O_X(K_X+H')) \leq 0.
    $$
    As $K_X.C = 0$ we have $K_X.H' = -2K_X.H$ and applying Riemann--Roch to the left hand side yields $H'^2-2K_X.H\leq 2q-2$. As $C$ is nef and $H'$ is effective, we have $C.H'= 2H.H'+H'^2\geq 0$. Together with the previous inequality we get $H.H'\geq 1-q-K_X.H$. Combining this with \eqref{eqn:lemmaproofeqn1} for $(\delta,n) = (1,3)$ we obtain
    $$
    d-l+K_X.H = (H^2+K_X.H) + H.H' \geq q-1 - H.K_X.
    $$
    On the other hand, for $(\delta,n) = (1,3)$, \eqref{eqn:lemmaproofeqn2} gives
    $$
    d-l+K_X.H \geq 2(2q-2)+\frac{1}{2}H.K_X.
    $$
    Multiplying the former bound by $\frac{1}{3}$ and the latter by $\frac{2}{3}$, their sum gives
    $$
    d-l+K_X.H\geq 3q-3\geq 3 \mathrm{gon}(B)-3
    $$
    since $q \geq 2$, as desired.
\end{proof}

\begin{lemma}\label{lem:contractionLemma}
    Let $X$ be a smooth projective ruled surface of irregularity $q > 0$. Denote by $f \colon X \rightarrow B$ the ruling map and let $F$ be the numerical class of a fiber. Assume the cycle of a reducible fiber of $f$ has the form $\Gamma_0 + \sum_i \mu_i \Gamma_i$, where $\mu_i \geq 1$ and $\Gamma_0$, $\Gamma_i$ are all smooth, irreducible rational curves of negative self-intersection, i.e., we assume that there is at least one component of multiplicity $1$, namely $\Gamma_0$. Then there is a $\Gamma_i$, $i \neq 0$, which is a $(-1)$-curve.
\end{lemma}

\begin{proof}
    Adjunction gives $K_X.F = -2$ because a general fiber of $f$ is isomorphic to $\P^1$ and $F^2 = 0$. Assume that $\Gamma_i^2 \leq -2$ for all $i \neq 0$. Adjunction implies $K_X.\Gamma_i = -2 - \Gamma_i^2 \geq 0$ and $K_X.\Gamma_0 = -2 - \Gamma_0^2 \geq -1$. Since $F \equiv \Gamma_0 + \sum_i \mu_i \Gamma_i$, we obtain
    \begin{equation*}
        -2 = K_X.F = K_X.\Gamma_0 + \sum_i \mu_i K_X.\Gamma_i \geq -1,
    \end{equation*}
    a contradiction.
\end{proof}

\begin{lemma}\label{lem:irregularityone}
    With the notation of Lemma~\ref{lem:DimensionofTinequality}, let $q = 1$ and let $f \colon X \to B$ be the ruling with $B$ an elliptic curve. Assume $n \geq 2$ and let $\mathcal{T} \subseteq \mathrm{Hilb}^l(X) \times \Pic(X)$ be an irreducible, locally closed stratum whose image $\mathcal{P} = \mathrm{im}(\mathcal{T} \to \Pic(X))$ has dimension $1$ and consists only of line bundles $H$ satisfying $H.(C-H+K_X) = c-2$. Then, $c = 2n-2$.
\end{lemma}

\begin{proof}
    The inequality $c \leq 2n-2$ holds automatically because $c \leq \gon(C) - 2 \leq n \gon(B) - 2 = 2n-2$ for $C \in |\mathcal{L}|_{\mathrm{sm}}$ general. For the sake of contradiction, we assume $c < 2n-2$.
    For a point $H \in \mathcal{P}$ of local dimension $1$, any local branch of $\mathcal{P}$ at $H$ is dominated by a family of the form $H_{\alpha} = H \otimes f^\ast(\alpha)$, where $\alpha \in \Pic^0(B)$ varies in a $1$-dimensional family over some irreducible base curve (that can be taken to be some open of $\Pic^0(B) \simeq B$ since $B$ is elliptic). We will establish the following properties for every such family $H_\alpha$:
    \begin{enumerate}
        \item $C - 2H_\alpha + K_X$ is effective for all $\alpha$;
        \item $H_\alpha$ has no fixed curves for all $\alpha$;
        \item $H_\alpha|_C$ computes $\Cliff(C)$ for all $\alpha$ and is hence basepoint-free;
        \item $h^0(H_\alpha|_C) = h^0(H_\alpha)$ for all $\alpha$.
    \end{enumerate}
    The second property holds because $H_\alpha$ appears in some Donagi--Morrison sequence. For the same reason, both $H_{\alpha}|_C$ and $(C-H_{\alpha})|_C$ contain a pencil. They are Serre-dual to each other since $C \cap E =\emptyset$, thus $H_{\alpha}|_C$ contributes to the Clifford index $\mathrm{Cliff}(C)$. For $C \in |\mathcal{L}|_{\mathrm{sm}}$ general, the Clifford index inequality for $H_{\alpha}|_C$ together with $H.(C-H+K_X) = c-2$ gives
    \begin{equation*}
        c \leq H.C - 2h^0(H_\alpha|_C) + 2 \leq H.(C-H+K_X) + 2 - 2h^1(H_\alpha) = c - 2h^1(H_\alpha).
    \end{equation*}
    Therefore, $h^1(X,\O_X(H_{\alpha})) = 0$ and $h^0(H_\alpha|_C) = h^0(H_\alpha)$ for all $\alpha$, and $H_{\alpha}|_C$ \emph{computes} $\mathrm{Cliff}(C)$. In particular, $H_{\alpha}|_C$ is basepoint-free for all $\alpha$, establishing properties $(3)$ and $(4)$.
    For property $(1)$, consider the short exact sequence
    \begin{equation*}
        0 \rightarrow \O_X(K_X-H_\alpha) \rightarrow \O_X(-H_\alpha) \rightarrow \O_E(-H_\alpha) \rightarrow 0.
    \end{equation*}
    Since $h^1(\O_X(K_X-H_\alpha)) = h^1(\O_X(H_\alpha)) = 0$ and $h^0(\O_X(-H_\alpha)) = 0$, we get $h^0(\O_E(-H_\alpha)) = 0$. Because $H_\alpha$ has no fixed curves, $\O_E(-2H_\alpha) \subseteq \O_E(-H_\alpha)$, hence $h^0(\O_E(-2H_\alpha)) = 0$. Letting $H'_\alpha \coloneqq C - 2H_\alpha \geq 0$, every global section of $\O_X(H_\alpha')$ vanishes along $E$ because $\O_X(C-2H_\alpha)|_E = \O_E(-2H_\alpha)$ since $C \cap E = \emptyset$ and $h^0(\O_E(-2H_\alpha)) = 0$ as we have just seen.

    In order to arrive at a contradiction, consider now all $1$-dimensional families $H_\alpha$ of the form $H \otimes f^\ast(\alpha)$ with $\alpha \in \Pic^0(B)$ varying in some open such that all $H_\alpha$ satisfy the four properties listed above. We do \emph{not} assume that every such $H_\alpha$ necessarily comes from a Donagi--Morrison sequence. We claim that the existence of such a family contradicts $c < 2n-2$. To see this, we may choose one such family $H_\alpha$ for which $N \coloneqq h^0(X,\O_X(H_\alpha))$ is generically minimal. After shrinking the family, we may assume $h^0(X, \O_X(H_\alpha))$ is constant.

    \noindent Setting $a\coloneqq H_{\alpha}.F$, we first claim that $N \leq a+2$. If not, then $N \geq a+3$ and we set $s \coloneqq \dim \mathrm{im}(H^0(F, \O_F(a)) \to H^1(X,\O_X(H_\alpha - F)))$ where, by abuse of notation, $F \simeq \P^1$ also denotes the divisor class of some general fiber. The divisor sequence
    $$
        0 \to \O_X(H_\alpha - F) \to \O_X(H_\alpha) \to \O_{\P^1}(a) \to 0
    $$
    yields $h^0(X,\O_X(H_\alpha - F)) = N-a-1+s \geq 2$ since $N \geq a+3$. Setting again $H'_\alpha \coloneqq C- 2H_\alpha \geq 0$, we also have $h^1(C,\O_X(H_\alpha - F)|_C) = h^0(C, \O_X(H_\alpha + H'_\alpha + F)|_C) \geq h^0(X, \O_X(H_\alpha + H'_\alpha + F)) \geq 2$, so $\O_X(H_\alpha-F)|_C$ contributes to $\mathrm{Cliff}(C)$. We get
    \begin{equation}\label{eqn:intermediatecliffineq}
        c \leq \mathrm{Cliff}(\O_X(H_\alpha-F)|_C) \leq H_\alpha.C - n - 2(N-a-2+s) = c-(n-2a-2)-2s,
    \end{equation}
    where the last equality follows from $c = \mathrm{Cliff}(\O_X(H_\alpha)|_C) = H_\alpha.C - 2(N-1)$ since $H_\alpha|_C$ computes $\mathrm{Cliff}(C)$. As $C-2H_\alpha+K_X$ is effective and $F$ is nef,
    \begin{equation}\label{eqn:intermediateeqn1}
        n-2a-2 = (C - 2H_\alpha + K_X).F \geq 0,
    \end{equation}
    and then \eqref{eqn:intermediatecliffineq} forces $n=2a+2$, $s=0$, $h^0(C, \O_X(H_\alpha-F)|_C) = h^0(X,\O_X(H_\alpha-F))$ and that $(H_\alpha - F)|_C$ \emph{computes} $\mathrm{Cliff}(C)$. Consider now the fixed divisor $P_\alpha$ of $|H_\alpha - F|$. We have $P_\alpha.C = 0$ as $(H_\alpha - F)|_C$ is basepoint-free (since it computes $\Cliff(C)$). Since $C$ is big and nef, there are only finitely many irreducible curves on $X$ intersecting $C$ trivially, and every component of $P_\alpha$ is one of them, for every $\alpha$. Hence, $P_\alpha$ is generically constant in $\alpha$. Denote its generic member by $P$. Then, shrinking the family if necessary, $\widetilde{H}_\alpha \coloneqq H_{\alpha}-F-P$ is still a $1$-dimensional family satisfying the four properties above. The only property that is not obvious for $\widetilde{H}_\alpha$ is $(4)$ which follows from $h^0(\widetilde{H}_\alpha) = h^0(H_\alpha - F) = h^0((H_\alpha-F)|_C) = h^0(\widetilde{H}_\alpha|_C)$. However, since $h^0(X,\O_X(\widetilde{H}_{\alpha})) = N-a-1 < N$, this contradicts the minimality of $N$. We therefore have $N \leq a+2$, as desired.

    \noindent For a family $H_\alpha$ satisfying the four properties above such that $h^0(X,\O_X(H_\alpha)) = N \leq a+2$ is minimal, we now claim $b \coloneqq H_\alpha.E \geq 2$. First of all, we have $E = \sigma_1 + \sigma_2 + D$ as divisors on $X$, where $D \geq 0$ is supported on fibers of $f \colon X \rightarrow B$ and $\sigma_1$, $\sigma_2$ are sections of $f$ (possibly equal). For this, consider the contraction $\pi \colon X \to S$ to a minimal, geometrically ruled surface $S$, factoring $f \colon X \xrightarrow{\pi} S \to B$. Then $\pi_\ast(E) \in |-K_S|$ because for the blow-up of a single point, the canonical class is the sum of the pullback of the old canonical class and the exceptional divisor. Now, it is well-known that every anticanonical divisor on $S$ satisfies the claimed decomposition. The section $\sigma_i \colon B \rightarrow S$ gives a rational section of $f$ which extends to a global section of $f$ because $B$ is smooth. Since $\pi$ is an isomorphism outside the exceptional divisors, this proves that $\sigma_1 + \sigma_2 \leq E$ as actual divisors on $X$ (not only as divisor \emph{classes}). Since $E.F = 2$, the divisor $D \coloneqq E-\sigma_1-\sigma_2 \geq 0$ satisfies $D.F = 0$, so $D$ is supported on fibers, as claimed. Now, $H_\alpha.\sigma_i \geq 0$ since $H_\alpha$ is nef. As $H_\alpha$ has no fixed curves, if $H_\alpha.\sigma_i = 0$, then $H_\alpha|_{\sigma_i} = \O_{\sigma_i} = \O_B$. Since $\sigma_i$ is a section of $f$ and we may write $H_\alpha = H \otimes f^\ast(\alpha)$, we have $\O_B = H_\alpha|_{\sigma_i} = H|_{\sigma_i} \otimes \alpha$ contradicting the variation of $\alpha$. Hence, $H_\alpha.\sigma_i \geq 1$ for $i = 1,2$ and therefore $b = H_\alpha.E \geq 2$, as desired.

    \noindent Next, let $a \coloneqq H_\alpha.F \geq 0$. We claim $H_\alpha^2 \geq ab$. Writing $H_\alpha = H \otimes f^\ast(\alpha)$ for $\alpha \in \Pic^0(B) \simeq B$ varying in some open, we consider again the contraction $\pi \colon X \to S$ and denote by $Z_i \subset X$ the (pullbacks to $X$ of the) exceptional divisors. For the successive proper transforms of $H_0 \coloneqq \pi_\ast(H)$ we write $H_i \coloneqq \pi_{i}^\ast(H_{i-1}) - m_i Z_i$. As the $Z_i$ are contained in fibers, we have $0 \leq m_i \leq a$. Indeed, intersecting $H_i$ with $Z_i$ gives $m_i$, and as $H_i$ is nef we have $m_i = H_i.Z_i \leq H_i.F_i = H.F = a$, where $F_i$ is the numerical class of the fiber containing $Z_i$ at the $i$-th blow-up step.
    Now, $H_0^2 = a(H_0.(-K_S))$ on $S$. Indeed, writing $H_0 \equiv a \sigma + h F_0$ for the numerical class, where $\sigma$ is the section class satisfying $\sigma^2 = -e$, we have $H_0^2 = a(2h-ea)$. Adjunction on the section and fiber give $K_S \equiv -2\sigma - eF_0$, so $H_0.(-K_S) = 2h-ea$, as desired. Next, under a single blow-up,
    \begin{equation*}
        H_i^2 - a(H_i.(-K_{X_i})) = H_{i-1}^2 - a(H_{i-1}.(-K_{X_{i-1}})) - m_i^2 + am_i.
    \end{equation*}
    Iteratively, using $b = H.E = H.(-K_X)$, we get
    \begin{equation}\label{eqn:H^2ab}
        H^2 - ab = \sum_i m_i(a-m_i) \geq 0,
    \end{equation}
    proving our claim. This shows that $2N = 2h^0(\O_X(H_\alpha)) \geq 2 \chi(\O_X(H_\alpha)) = H_\alpha^2 + H_\alpha.(-K_X) \geq ab + b = (a+1)b$, i.e., $N \geq \frac{b}{2}(a+1)$.

    \noindent Next, we show $a \geq 1$. Indeed, if $a = 0$, then each $H_{\alpha} = f^\ast(\beta)$ is a pullback from $B$. Since $f_\ast \O_X = \O_B$, we have $f_\ast(\O_X(H_\alpha)) = \beta$, in particular
    \begin{equation*}
        h^0(B,\beta) = h^0(X, \O_X(H_\alpha)) = h^0(C, \O_X(H_\alpha)|_C) \geq 2
    \end{equation*}
    as $H_{\alpha}|_C$ computes $\mathrm{Cliff}(C)$ for $C$ general. In particular, since $B$ is elliptic, $\deg(\beta) \geq 2$. Moreover, the Clifford equality for $H_\alpha|_C$ yields
    \begin{equation*}
        c = C.H_\alpha - 2h^0(H_\alpha|_C) + 2 = \deg(\beta) n - 2\chi(H_\alpha) + 2 = \deg(\beta) n - 2\deg(\beta) + 2 \geq 2n-2,
    \end{equation*}
    using $n \geq 2$ and Riemann--Roch for $\chi(H_\alpha)$ together with $H_\alpha = f^\ast(\beta) \equiv \deg(\beta)F$ and $F.(-K_X) = 2$. This contradicts $c < 2n-2$.

    \medskip
    \noindent The inequalities $b \geq 2$, $a \geq 1$, $N \leq a+2$ and $N \geq \frac{b}{2}(a+1)$ only leave the three possibilities
    \begin{equation}
        (b,N) \in \{(2,a+1), (3,3), (2,a+2)\}.
    \end{equation}
    Moreover, in the second case necessarily $a = 1$, and in the first two cases one has equality in $N \geq \frac{b}{2}(a+1)$. This implies that $h^0(\O_X(H_\alpha)) = \chi(\O_X(H_\alpha))$, i.e., $h^1(\O_X(H_\alpha)) = 0$, and that $\sum_i m_i (a-m_i) = 0$, so that all $m_i \in \{0,a\}$. In the third case, expressing $H^2$ through Riemann--Roch, necessarily
    \begin{equation*}
        2 = 2h^1(\O_X(H_\alpha)) + \sum_i m_i(a-m_i),
    \end{equation*}
    so either $h^1(\O_X(H_\alpha)) = 1$ and all $m_i \in \{0,a\}$ or $h^1(\O_X(H_\alpha)) = 0$ and at least one $m_i$ lies strictly between $0$ and $a$. This implies $a \geq 2$, and since $m_i(a-m_i) \geq a-1$ whenever $1 \leq m_i \leq a-1$, necessarily $a-1 \leq 2$, so $a \in \{2,3\}$.

    \noindent Before proceeding to the case distinction, we make three observations:
    \begin{enumerate}[(i)]
        \item Writing our family as $H_\alpha = H \otimes f^\ast(\alpha)$ for $\alpha \in \Pic^0(B) \simeq B$ varying in some open, the properties $(1)$ and $(4)$ actually extend to all of $\Pic^0(B) \simeq B$. This follows from upper semicontinuity of cohomology: We observe that $h^0(\O_X(H_\alpha)) \geq N$ for \emph{all} $\alpha \in \Pic^0(B)$, so that $H_\alpha|_C$ contributes to Clifford for \emph{every} $\alpha \in \Pic^0(B)$. For $C$ general, we get
        \begin{equation*}
            c \leq \Cliff(H_\alpha|_C)  = H_\alpha.C - 2h^0(H_\alpha|_C) + 2 \leq H_\alpha.C - 2 h^0(H_\alpha) + 2 \leq H_\alpha.C - 2(N-1) = c,
        \end{equation*}
        hence necessarily $h^0(\O_X(H_\alpha)|_C) =  h^0(\O_X(H_\alpha)) = N$ and $H_\alpha|_C$ computes the Clifford index and is thus basepoint-free for \emph{all} $\alpha \in \Pic^0(B)$.
        \item Assume $h^1(\O_X(H_\alpha)) = 0$ and let $F$ be the divisor class of a fiber of $f$. Then $h^2(\O_X(H_\alpha - F)) = h^0(\O_X(K_X + F - H_\alpha)) = 0$ as $F$ is nef and $F.(K_X + F - H_\alpha) = -2 + 0 - a < 0$. In particular, $h^0(\O_X(H_\alpha - F)) = N - a -1 + h^1(\O_X(H_\alpha - F))$. If the latter is $\geq 2$, then $(H_\alpha - F)|_C$ contributes to Clifford and so
        \begin{equation*}
            c \leq C.H_\alpha - n - 2(N-a-2) - 2h^1(\O_X(H_\alpha - F)) = c - (n - 2a -2) - 2h^1(\O_X(H_\alpha - F)).
        \end{equation*}
        So we have shown that, if $h^1(\O_X(H_\alpha)) = 0$ and $h^0(\O_X(H_\alpha - F)) \geq 2$, then $n = 2a+2$ and $h^1(\O_X(H_\alpha - F)) = 0$.
        \item If all $m_i \in \{0,a\}$, then we claim that there is a blow-down map $\rho \colon X \rightarrow Y$ to a geometrically ruled surface $Y$ over $B$ such that $H_\alpha = \rho^\ast(H_\alpha^Y)$ for $H_\alpha^Y = \rho_\ast(H_\alpha) \in \Pic(Y)$. Moreover, $-K_Y$ is effective since $\rho_\ast(E) \in |-K_Y|$. Letting $\sigma$ denote the class of a minimal section of the ruling $Y \rightarrow B$ and $F_0$ the fiber class, we have $\sigma^2 = -e$ for some $e \geq 0$.
        
        \noindent To see this, we first observe that $H_\alpha = H \otimes f^\ast(\alpha)$ intersects all but one irreducible component of every fiber of $f$ trivially, and for the remaining component $\Gamma$ we have $H_\alpha.\Gamma = a$. Moreover, $\Gamma$ has multiplicity $1$ in its fiber. Indeed, these claims hold for $H_0$ on the minimal surface $S$ from above, so we may proceed by induction. If the claim holds before the $i$-th blow-up $\pi_i$, assume first that $m_i = 0$. In this case, $H_i = \pi_i^\ast(H_{i-1})$ and $H_i.Z_i = 0$. In particular, if $Z \neq Z_i$ is any irreducible component of some fiber, then $H_i.\widetilde{Z} = H_i.\pi^\ast(Z) = H_{i-1}.Z$, so the claim is still true after the $i$-th blow-up. If instead $m_i = a$, then $H_i.Z_i = a$, $Z_i$ has multiplicity $1$ in its fiber and necessarily $H_i$ intersects every other component of the same fiber trivially because $H_i.F_i = H.F = a$ and $H_i$ is nef. This proves the claim and puts us into the situation of Lemma~\ref{lem:contractionLemma}. We can therefore contract a $(-1)$-curve contained in some fiber of $f$ that intersects $H_\alpha$ trivially. By the above argument in the case $m_i = 0$, this contraction does not change the intersection of $H_\alpha$ with the other components of this or any other fiber of $f$. The argument may therefore be repeated until every fiber is irreducible, i.e., until we arrive at a geometrically ruled surface $Y$ (which may be different from $S$). We still need to justify that $e \geq 0$. In general, $e \geq -1$, so we only need to exclude the case $e = -1$. Indeed, $-K_Y \equiv 2 \sigma - F_0$. On the other hand, we have seen above that $\rho_\ast(E) = \sigma_1 + \sigma_2 + D$ for some $D \geq 0$ supported on fibers and two sections $\sigma_i$ of $Y \rightarrow B$. Write $D \equiv r F_0$ for some $r \geq 0$. As $\sigma$ is a minimal section, $\sigma_i^2 \geq \sigma^2 = -e = 1$. Writing $\sigma_i \equiv \sigma + s_i F_0$, this shows $s_i \geq 0$ and hence
    \begin{equation*}
        2 \sigma - F_0 \equiv -K_Y \sim \rho_\ast(E) = \sigma_1 + \sigma_2 + D \equiv 2 \sigma + (s_1 + s_2 + r) F_0,
    \end{equation*}
    a contradiction. Hence, $e \geq 0$.
    \end{enumerate}
    
    \noindent \textbf{Case} $(b,N) = (2,a+1)$. Since all $m_i \in \{0,a\}$, we may apply observation (iii). Write $H_\alpha^Y \equiv a \sigma + (ae+t)F_0$, where $t = H_\alpha^Y.\sigma \geq 0$ since $H_\alpha^Y$ is nef. In fact $t \geq 1$ as otherwise $t = H_\alpha^Y.\sigma = 0$, so also $H_\alpha.\sigma = 0$ on $X$ using that the section $\sigma$ of $Y \rightarrow B$ extends to a section on $X$ of $f$. Since $H_\alpha$ has no fixed curves, this implies $\O_B = H_\alpha|_\sigma = H|_\sigma \otimes \alpha$, contradicting the variation of $\alpha$. Since $-K_Y \equiv 2 \sigma + e F_0$, we get $2 = b = H_\alpha^Y.(-K_Y) = ae+2t$. Since $e \geq 0$ and $t \geq 1$, the only possible case is $e = 0$ and $t = 1$. In particular, $\sigma^2 = 0$ and $H_\alpha^Y.\sigma = 1$. Let $C' \coloneqq \rho(C)$. Write $C' \equiv n \sigma + r F_0$. Then $2nr = C'^2 \geq C^2 > 0$, so $r \geq 1$. In particular $C'.\sigma = r \geq 1$, i.e., $C' \cap \sigma \neq \emptyset$. Let $p' \in C' \cap \sigma$ and $p \in  C$ a point over $p'$. Again using $\sigma$ both for the section of $Y$ and its extension to $X$, we have $H_\alpha^Y|_\sigma = H_\alpha|_\sigma = H|_\sigma \otimes \alpha$. This is a line bundle of degree $1$ on $B$ and is therefore of the form $\O_B(q)$ for some $q$ (which depends on $\alpha$). Fixing some $\alpha$ and hence $q$, twisting by $\O_B(p-q) \in \Pic^0(B)$, we obtain the line bundle $H_{\alpha'}$ on $X$ with $\alpha' \in \Pic^0(B)$ potentially lying outside of the family. Nonetheless, observation (i) shows that $H_{\alpha'}|_C$ still computes $\Cliff(C)$ for $C$ general and is hence basepoint-free. On the other hand, $H_{\alpha'}^Y|_\sigma = \O_B(p)$, so that all global sections of $\O_Y(H^Y_{\alpha'})$ vanish at $p \in Y$. We have $H^0(\O_X(H_{\alpha'})) = H^0(\O_Y(H^Y_{\alpha'}))$ since $\rho_\ast(\O_X) = \O_Y$ which implies that every global section in $H^0(\O_X(H_{\alpha'}))$ vanishes at every preimage in $X$ of $p \in Y$.  In particular, since $h^0(\O_X(H_{\alpha'})|_C) = h^0(\O_X(H_{\alpha'}))$ by observation (i), every preimage of $p$ in $C$, which exists because $p \in \rho(C)$, is a basepoint of $H_{\alpha'}|_C$, a contradiction.

    \noindent \textbf{Case} $(b,N) = (2,a+2)$. If $h^1(\O_X(H_\alpha)) = 1$, then again all $m_i \in \{0,a\}$ and we may repeat the argument from the previous case obtaining the same contradiction. Otherwise, $h^1(\O_X(H_\alpha)) = 0$ and $\sum_i m_i (a-m_i) = 2$. Here, we subtract a general fiber $F$ from $H_\alpha$. As in observation (ii), $h^2(\O_X(H_\alpha - F)) = 0$, so $h^0(\O_X(H_\alpha - F)) = (N-a-1) + h^1(\O_X(H_\alpha - F)) = 1 + h^1(\O_X(H_\alpha - F))$. Hence, if $h^0(\O_X(H_\alpha - F)) \geq 2$, then $h^1(\O_X(H_\alpha - F)) \geq 1$, contradicting observation (ii). Therefore, $h^0(\O_X(H_\alpha - F)) = 1$ and $h^1(\O_X(H_\alpha - F)) = 0$. Moreover, this must in fact hold for \emph{all} $\alpha \in \Pic^0(B)$ by combining the previous argument with observations (i) and (ii). Now, $h^2(\O_X(H_\alpha - F - E)) = h^0(\O_X(F - H_\alpha)) = 0$. The short exact sequence
    \begin{equation*}
        0 \rightarrow \O_X(H_\alpha - F - E) \rightarrow \O_X(H_\alpha - F) \rightarrow \O_X(H_\alpha - F)|_E \rightarrow 0
    \end{equation*}
    then shows that $h^1(E,\O_X(H_\alpha - F)|_E) = 0$, again for \emph{every} $\alpha \in \Pic^0(B)$. Since $(H_\alpha - F).E = b - 2 = 0$ and $E$ is Gorenstein of arithmetic genus $1$, this implies $h^0(E,\O_X(H_\alpha - F)|_E) = 0$ by Riemann--Roch. Consider $V_i \coloneqq R^i(f|_E)_\ast(\O_X(H_\alpha - F)|_E)$ for $i \in \{0,1\}$. Then $V_1$ is torsion on $B$ because $f|_E$ is generically finite. Since $H^0(V_1) \simeq H^1(E, \O_X(H_\alpha - F)|_E) = 0$ by the Leray spectral sequence, $V_1 = 0$. Moreover, $V_0$ has rank $2$ because $\deg(f|_E) = E.F = 2$. Since $h^0(V_0) = h^0(E,\O_X(H_\alpha - F)|_E) = 0$, $V_0$ is torsion-free and hence a vector bundle. Again by the Leray spectral sequence, $H^i(V_0) = H^i(E, \O_X(H_\alpha - F)|_E) = 0$ for $i = 0,1$. Hence, choosing any point $p \in B$, the short exact sequence
    \begin{equation*}
        0 \rightarrow V_0 \rightarrow V_0(p) \rightarrow V_0(p)|_p \rightarrow 0
    \end{equation*}
    shows that $H^0(V_0(p)) \rightarrow H^0(V_0(p)|_p) \simeq \C^2$ is an isomorphism, in particular $h^0(V_0(p)) = 2$. The evaluation map $H^0(V_0(p)) \rightarrow H^0(V_0(p)|_q) \simeq \C^2$ is in fact an isomorphism for every $q \in B$ because its kernel is $H^0(V_0(p-q))$, and $V_0(p-q) = (f|_E)_\ast(\O_X(H_\alpha - F)\otimes f^\ast \O_B(p-q)|_E)$ has no global sections because $\O_B(p-q)$ can be absorbed into $\alpha$. This implies that $V_0(p)$ is globally generated, i.e., $\O_B^2 \rightarrow V_0(p)$ is surjective and hence an isomorphism. In particular, $\deg(V_0(p)) = 0$, contradicting $\deg(V_0(p)) = \deg(V_0) + 2 = \chi(B,V_0) + 2 = \chi(E,\O_X(H_\alpha - F)|_E) + 2 = 2$.

    \noindent \textbf{Case} $(a,b,N) = (1,3,3)$. Here, again $h^1(\O_X(H_\alpha)) = 0$ and all $m_i \in \{0,a\}$, so we may use observation (iii). With the same notation as in the case $(b,N) = (2,a+1)$, we have $H_\alpha^Y \equiv a \sigma + (ae+t) F_0$, and again $t \geq 1$ by the same argument. Here, we get $3 = b = H_\alpha^Y.(-K_Y) = e+2t$ since $a = 1$ and therefore $e = t = 1$, so $\sigma^2 = -1$ and $H_\alpha^Y \equiv \sigma + 2F_0$. Then $H_\alpha^Y.\sigma = 1$ and, again extending $\sigma$ to a section of $f$, also $H_\alpha.\sigma = 1$. For $\alpha$ sufficiently general, the degree $1$ line bundle $H_\alpha|_\sigma$ is isomorphic to $\O_B(p)$, where $p$ is such that the fiber $F_p \coloneqq f^{-1}(p)$ over $p$ is isomorphic to $\P^1$. If not, then $p$ could only be one of finitely many points, contradicting the fact that $H_\alpha|_\sigma = H|_\sigma \otimes \alpha$ varies in a $1$-dimensional family. Hence, all global sections of $H_\alpha$ vanish at the point $P \in \sigma \subseteq X$ corresponding to $p \in B$ and so the restriction map $H^0(X,\O_X(H_\alpha)) \rightarrow H^0(F_p,\O_{\P^1}(1))$ lands inside the $1$-dimensional subspace of sections vanishing at $P \in F_p$. In particular, since $h^0(\O_X(H_\alpha)) = N = 3$, this shows that $h^0(\O_X(H_\alpha - F_p)) \geq 2$. From observation (ii) we obtain $n = 2a+2 = 4$ and $h^1(\O_X(H_\alpha - F_p)) = 0$. Let $r \coloneqq h^0(\O_X(H_\alpha)) - h^0(\O_X(H_\alpha - F_p)) \in \{0,1\}$. Since $(H_\alpha - F_p)|_C$ contributes to the Clifford index, for $C$ general, we have
    \begin{align*}
        c &\leq C.H_\alpha - n - 2(h^0(\O_X(H_\alpha)) - r) + 2 \\
        &= C.H_\alpha - 2(N-1) - n + 2r \\
        &= c - n + 2r \leq c - 2,
    \end{align*}
    a contradiction. Since there are no cases left, this concludes the proof.
\end{proof}

\begin{proposition}\label{prop:genericgreen2}
    Keeping the assumptions of this subsection, let $\mathcal{W} \subseteq \widetilde{\mathcal{W}}^{1,\mathrm{bpf}}_d(|\mathcal{L}|_{\mathrm{sm}})$ be a dominating irreducible component such that, for the general $(C,A) \in \mathcal{W}$, the Lazarsfeld--Mukai bundle $F_{C,A}$ is non-simple but indecomposable. If $q \geq 1$, let $f \colon X \to B$ be the ruling and assume $n \geq 2$. Then, the general fiber of $\mathcal{W}$ over $|\mathcal{L}|_{\mathrm{sm}}$ satisfies the linear growth condition \eqref{eqn:lineargrowthcondition} if one of the following conditions holds:
    \begin{enumerate}[(i)]
        \item $q = 0$;
        \item $q \geq 1$ and $n  \gon(B) \neq c + 2$;
        \item $q \geq 2$, $n  \gon(B) = c + 2$ and $B$ has only finitely many gonal pencils.
    \end{enumerate}
\end{proposition}

\begin{proof}
    Let $(C,A)\in \mathcal{W}$ be general. As $F = F_{C,A}$ is non-simple but indecomposable, Proposition~\ref{prop:DMextension} yields that $F$ sits in a Donagi--Morrison extension
    \begin{equation}\label{eqn:DMExtension}
    0 \to \O_X(H+H') \to F \to \mathcal{I}_Z(H) \to 0,
    \end{equation}
    where $H > 0$, $H' \geq 0$, $C \sim 2H+H'$ and $Z \subset X$ is a $0$-dimensional lci subscheme. We bound the dimension of the set of pairs $(F,V)$, where $F$ is the isomorphism class of a Lazarsfeld--Mukai bundle sitting in such an extension and $V \subseteq H^0(X,F)$ is a $2$-dimensional subspace. Restricting \eqref{eqn:DMExtension} to $E$ and using $C \cap E = \emptyset$ as well as the torsion-freeness of $\mathcal{I}_Z(H)$, we obtain
    $$
    0 \to \O_E(-H) \to \O_E^2 \to \mathcal{I}_Z(H)|_E \to 0.
    $$
    As $E$ is Gorenstein of arithmetic genus~$1$, we have $\omega_E \simeq \O_E$, so Riemann--Roch implies that any rank~$1$ coherent sheaf $\mathcal{F}$ on $E$ satisfies $\chi(\mathcal{F}) = \deg(\mathcal{F})$. Hence,
    \begin{equation*}
        \chi(\mathcal{I}_Z(H)|_E,\O_E(-H)) = -\chi(\mathcal{I}_Z(2H)|_E) = -\deg(\mathcal{I}_Z(2H)|_E) =  2H.K_X,
    \end{equation*}
    and thus
    \begin{equation*}
    \ext^1(\mathcal{I}_Z(H)|_E,\O_E(-H)) = \hom(\mathcal{I}_Z(H)|_E,\O_E(-H)) - 2H.K_X.
    \end{equation*}
    We want to bound the dimension of the space of extension classes in $\Ext^1(\mathcal{I}_Z(H)|_E,\O_E(-H))$ whose middle term is isomorphic to $\O_E^2$. We do this by bounding the dimension of the space of actual extension sequences, parametrized by pairs $(\alpha,\varphi)$, where $\alpha\colon \O_E^2\twoheadrightarrow \mathcal{I}_Z(H)|_E$ is a surjection and $\varphi\colon \ker(\alpha)\xrightarrow{\sim}\O_E(-H)$ an isomorphism, and then passing from a sequence to its equivalence class in $\Ext^1(\mathcal{I}_Z(H)|_E,\O_E(-H))$, thereby subtracting from the count of pairs $(\alpha, \varphi)$ the dimension of the fibers, i.e., the dimension of $\Aut(\O_E^2)/\mathrm{Stab}_{(\alpha,\varphi)}$.
    
    \noindent First of all, using again $\chi(\mathcal{I}_Z(H)|_E) = \mathrm{deg}(\mathcal{I}_Z(H)|_E) = -H.K_X$, we find $h^0(E,\mathcal{I}_Z(H)|_E) = -H.K_X+h^1(E,\mathcal{I}_Z(H)|_E)$. Thus, the count of $\alpha$ contributes at most $-2H.K_X+2h^1(E,\mathcal{I}_Z(H)|_E)$ moduli. As $H|_E$  is a line bundle on each irreducible component of $E$, we have $\mathrm{Aut}(\O_X(H)|_E) \simeq H^0(E,\O_E)^{\times}$, so $\varphi$ contributes at most $h^0(E,\O_E) = 1 + q$\footnote{In the irregularity zero case, the equality $h^0(E,\O_E) = 1 + q$ is clear. If $q > 0$, then this follows from the ruling on $X$ as in the proof of Proposition \ref{prop:genericgreen1}.} moduli. Now, $\tau \in \mathrm{Aut}(\O_E^2)$ acts on the parameter space of pairs $(\alpha,\varphi)$ by $\tau \cdot (\alpha, \varphi) = (\alpha \circ \tau, \varphi \circ \tau)$.
    By a general fact of homological algebra, the stabilizer in $\Aut(\O_E^2)$ of the extension short exact sequence associated to $(\alpha, \varphi)$ is isomorphic to $\Hom(\mathcal{I}_Z(H)|_E,\O_E(-H))$. Hence, setting
    \begin{equation*}
        \mathcal{V}_E \coloneqq \{[0 \to \O_E(-H) \to M \to \mathcal{I}_Z(H)|_E \to 0] \in \Ext^1(\mathcal{I}_Z(H)|_E,\O_E(-H)) \  | \ M \simeq \O_E^2\},
    \end{equation*}
    we obtain
    \begin{align*}
    \dim \mathcal{V}_E &\leq 2h^0(E,\mathcal{I}_Z(H)|_E) + \mathrm{aut}(\O_E(-H)) - \big( \mathrm{aut}(\O_E^2)- \hom(\mathcal{I}_Z(H)|_E,\O_E(-H)) \big) \\
    &= -2H.K_X+2h^1(E,\mathcal{I}_Z(H)|_E)-3(h^1(X,\O_X)+1)+\hom(\mathcal{I}_Z(H)|_E,\O_E(-H)).
    \end{align*}
    Therefore,
    \begin{equation}\label{eqn:codimExtensions}
        \gamma \coloneqq \mathrm{codim}(\mathcal{V}_E, \Ext^1(\mathcal{I}_Z(H)|_E,\O_E(-H))) \geq 3(q+1)-2h^1(E,\mathcal{I}_Z(H)|_E).
    \end{equation}
    Consider the restriction map
    \begin{equation}\label{eqn:restrictionofextensions}
    \rho \colon \Ext^1_X(\mathcal{I}_Z(H),\O_X(H+H')) \to \Ext^1_E(\mathcal{I}_Z(H)|_E,\O_E(-H)).
    \end{equation}
    By derived push-pull-adjunction and because $\mathcal{T}or_1^{\O_X}(\mathcal{I}_Z(H), \O_E) = \ker(\mathcal{I}_Z(-E) \rightarrow \mathcal{I}_Z) = 0$, we have 
    $$    \Ext^1_X(\mathcal{I}_Z(H),i_*\O_X(H+H')|_E)\simeq \Ext^1_E(\mathcal{I}_Z(H)|_E,\O_E(-H))
    $$
    using that $(H+H')|_E= -H|_E$ since $C \sim 2H + H'$ and $C \cap E = \emptyset$. Applying $\Hom(\mathcal{I}_Z(H),-)$ to the restriction sequence
    $$
    0 \to \O_X(H+H'+K_X) \to \O_X(H+H') \to \O_E(-H) \to 0
    $$
    shows that \eqref{eqn:restrictionofextensions} is surjective because 
    $$
    \Ext^2(\mathcal{I}_Z(H),\O_X(H+H'+K_X)) \simeq \Hom(\O_X(H+H'),\mathcal{I}_Z(H))^\vee = 0
    $$
    if $H' \neq 0$ or $Z \neq \emptyset$. If, on the other hand, $H' = 0$ and $Z = \emptyset$, then $\Ext^1(\O_X(H),\O_X(H)) \simeq H^1(X,\O_X)$ and $\Ext^1(\O_E(H),\O_E(H)) \simeq H^1(E,\O_E)$, and $\rho$ is the canonical map $H^1(X, \O_X) \rightarrow H^1(E, \O_E)$ which is injective. Now, if the middle term of an extension restricts to $F|_E \simeq \O_E^2$ on $E$, then one can check that the restricted sequence must split, hence so does the Donagi--Morrison sequence, i.e., the Lazarsfeld--Mukai bundle $F$ is decomposable, a contradiction. Therefore, for fixed $Z$ and $H$, we find that the space
    $$
    \mathcal{V}_{Z,H} \coloneqq \{[0 \to \O_X(H+H') \to F \to \mathcal{I}_Z(H) \to 0]\in \Ext^1(\mathcal{I}_Z(H),\O_X(H+H')) \ | \ F|_E\simeq \O_E^2\}
    $$
    of vector bundles in a Donagi--Morrison extension which become trivial on $E$ has codimension at least $\gamma \geq 3(q+1)-2h^1(E,\mathcal{I}_Z(H)|_E)$ in $\Ext^1(\mathcal{I}_Z(H),\O_X(H+H'))$.

    \noindent Next, we count parameters for $(Z,H)$. As in Lemma~\ref{lem:DimensionofTinequality}, consider an irreducible, locally closed stratum $\mathcal{T} \subseteq \mathrm{Hilb}^l(X) \times \Pic(X)$ and its image $\mathcal{P} = \mathrm{im}(\mathcal{T} \rightarrow \Pic(X))$. Set $p \coloneqq \dim \mathcal{P}$. Then, by Lemma~\ref{lem:DimensionofTinequality},
    \begin{equation}\label{eqn:boundDimT}
        \dim \mathcal{T} \leq p + 2l + 2(h^1(E,\O_E(H)) - h^1(E,\mathcal{I}_Z(H)|_E)).
    \end{equation}
    Recall that $l = d - H.(H+H')$. We may choose $\mathcal{T}$ such that the general element in $\mathcal{W}$ has an associated Lazarsfeld--Mukai bundle which sits in a Donagi--Morrison sequence for which $(Z,H)$ lies in $\mathcal{T}$. Combining \eqref{eqn:boundDimT} with the parameter count for extensions with fixed $(Z,H)$ as well as with the Grassmannian dimension for $2$-dimensional subspaces of global sections via Lemma~\ref{lem:constancyofLMsections}, we obtain the dimension estimate
    \begin{align}
    \begin{split}\label{eqn:T_jestimate}
        \dim \mathcal{W} &\leq \dim \mathcal{V}_{Z,H} + \dim \Gr(2,g-d+3) - (1 + h^0(X,\O_X(H'))) + \dim \mathcal{T} \\
        &\leq \ext^1(\mathcal{I}_Z(H),\O_X(H+H')) - \gamma + 2(g-d+1) - (1 + h^0(X,\O_X(H'))) + \dim \mathcal{T}\\
        &\leq g + (d-c-2) + p - (\eta-c)\\
        &= \dim |\mathcal{L}| + (d-c-2) + p - (\eta-c),
    \end{split}
    \end{align}
    where $\eta \coloneqq d-l + K_X.H + 2(q - h^1(E,\O_E(H)))$. Here, we subtract $1 + h^0(X,\O_X(H')) = \dim \Aut(F)$ because this is a lower bound on the dimension of every fiber of the map sending a pair $(F,V)$ with $V \subseteq H^0(F)$ a $2$-dimensional subspace to the quotient in the dual of the evaluation sequence
    \begin{equation*}
        0 \to F^\vee \to V^\vee \otimes \O_X \to i_\ast A \to 0,
    \end{equation*}
    where $i \colon C \hookrightarrow X$ is the inclusion of a curve in $|\mathcal{L}|_{\mathrm{sm}}$, i.e., $(F,V) \mapsto (C,A) \in \mathcal{W}$.
    
    \noindent It is left to prove that $p-(\eta-c)\leq 0$, as then $(\ref{eqn:T_jestimate})$ establishes the linear growth condition~\eqref{eqn:lineargrowthcondition}. For this, projecting the pencil $H^0(C,A)^\vee \subseteq H^0(F)$ to $H^0(\O_X(H))$ shows $h^0(C,\O_C(H)) \geq 2$. As also $h^1(C,\O_C(H)) = h^0(C,\O_C(H+H')) \geq h^0(X,\O_X(H+H')) \geq h^0(X,\O_X(H)) \geq 2$, we find that $\O_C(H)$ contributes to $\mathrm{Cliff}(C)$. Using Riemann--Roch on $X$, the Clifford inequality yields 
    \begin{align*}
    c \leq \mathrm{Cliff}(\O_C(H)) &\leq H.(H+H') + K_X.H + 2(q - h^1(X,\O_X(H)))\\
    &= d-l + K_X.H + 2(q - h^1(X,\O_X(H))).
    \end{align*}
    Restriction gives a surjection $H^1(X,\O_X(H))\twoheadrightarrow H^1(E,\O_E(H))$ and hence we have shown $c\leq \eta$. Thus, the inequality (\ref{eqn:T_jestimate}) already establishes linear growth (\ref{eqn:lineargrowthcondition}) if $p = 0$. This proves $(i)$.
    Moreover, if on $\mathcal{P}$ we have $-K_X.H = 0$, then $\O_X(H)|_E \simeq \O_E$ since $H$ has no fixed curve. As $\Pic(X) \to \Pic(E)$ has $0$-dimensional kernel because the tangent map $H^1(X,\O_X) \to H^1(E,\O_E)$ is injective, this implies $p = 0$ and the previous argument concludes the proof.
    We can therefore assume $-K_X.H > 0$ for $H \in \mathcal{P}$ and that $p \geq 1$.
    If $q = 1$, then by Lemma~\ref{lem:PicardBound} the only case for which \eqref{eqn:T_jestimate} does not establish linear growth is $p = 1$, $h^1(E,\O_E(H)) = 0$ as well as $\eta = c$ or, equivalently, $h^1(X,\O_X(H)) = h^1(E,\O_E(H)) = 0$ and $H.(C+K_X-H) = c-2$. This puts us into the situation of Lemma~\ref{lem:irregularityone} and we obtain $c = 2n-2$ with $n = \mathrm{deg}(f|_C)$, proving the $q = 1$ case of $(ii)$.

    \noindent Let now $q \geq 2$. By Lemma~\ref{lem:PicardBound}, we have $p + h^1(E,\O_E(H)) \leq q$. Setting
    \begin{equation*}
        \nu \coloneqq d-l + K_X.H
    \end{equation*}
    and using \eqref{eqn:T_jestimate}, it suffices to show $\nu \geq c-1$ as then $\eta - c = (\eta - \nu) + (\nu - c) \geq 2p - 1$ and so $p - (\eta - c) \leq 1-p \leq 0$, as desired.
    
    \noindent Let $f \colon X \to B$ be the ruling, $F$ the class of a fiber and let again $n = F.C = \deg(f|_C)$. We first assume $H.F > 0$ for $H \in \mathcal{P}$. We always have $c \leq \gon(C) - 2 \leq n \gon(B) - 2$, and so Lemma~\ref{lem:pullbackreduction} implies $\nu \geq c-1$.
    
    \noindent The only case left is $H.F = 0$. Then every component of $H$ is contained in a fiber, and the part $H_b$ of $H$ which is supported on a given fiber $F_b \coloneqq f^{-1}(b)$ intersects every component of $F_b$ trivially since $H$ and $F$ are both nef. In particular, $H_b^2 = 0$, so Zariski's lemma on the intersection matrix of the components of a fiber implies that $H = f^\ast L$ for some $L \in \Pic(B)$. In this case, the proof of \cite[Proposition 5.2]{GeigerKretschmerBN} shows that $Z$ is a subscheme of a $0$-dimensional complete intersection $Z'$ given by the vanishing of two sections of $\O_X(H)$. Since $H^2 = 0$, we have $Z = Z' = \emptyset$, so $l = 0$. The same proof then shows $A = H|_C$ and $h^0(X, \mathcal{L} \otimes f^\ast(L)^{-2}) > 0$. Hence, writing $t \coloneqq \deg L$ and $d = \deg A = nt$, if the condition
    \begin{equation}\label{eqn:gonalbasepencils}
    \dim \L\{L\in \Pic^t(B)\colon L \text{ is a complete, bpf pencil with } h^0(X,\mathcal{L}\otimes f^*L^{-2})>0\R\}\leq nt-c-2
    \end{equation}
    is satisfied whenever $nt \leq g-c$, then the linear growth condition \eqref{eqn:lineargrowthcondition} holds for the general fiber of $\mathcal{W}$ over $|\mathcal{L}|$. By \cite{FultonHarrisLazarsfeld}, the bound $\dim W^1_{\gon(B)}(B) \leq 1$ always holds, so \eqref{eqn:gonalbasepencils} is actually automatic when $n \gon(B) \neq c+2$. This is because $c+2 \leq \gon(C) \leq n \gon(B)$, whence, if equality does not hold, $n \gon(B) - c - 2 \geq 1$ follows. In the case of equality $n \gon(B) = c+2$, since $n \geq 2$, \cite{FultonHarrisLazarsfeld} still shows that \eqref{eqn:gonalbasepencils} is implied by $B$ having only finitely many gonal pencils.
\end{proof}

\begin{proposition}\label{prop:genericgreen3}
    Let $\mathcal{W} \subseteq \widetilde{\mathcal{W}}^{1,\mathrm{bpf}}_d(|\mathcal{L}|_{\mathrm{sm}})$ be an irreducible component such that, for the general $(C,A) \in \mathcal{W}$, the Lazarsfeld--Mukai bundle $F_{C,A}$ is a direct sum of two line bundles. Then the general fiber of $\mathcal{W}$ over $|\mathcal{L}|_{\mathrm{sm}}$ satisfies the linear growth condition \eqref{eqn:lineargrowthcondition}.
\end{proposition}

\begin{proof}
    Let $F = F_{C,A} = M \oplus N$ be a direct sum of line bundles. Then necessarily both $M$ and $N$ are generated away from finitely many points and hence nef. As $F|_E \simeq \O_E^2$ we have $(M+N).E=0$ and therefore $M|_E \simeq N|_E \simeq \O_E$. Again, as $\Pic(X) \to \Pic(E)$ has $0$-dimensional kernel, the isomorphism classes of $M$ and $N$ are locally constant. We first assume $M \neq N$. Up to exchanging $M$ and $N$ we may therefore assume $h^0(X,M^{\vee} \otimes N) = 0$. The short exact sequence
    $$
    0 \to \omega_X \otimes M \otimes N^{\vee} \to M \otimes N^{\vee} \to \O_E \to 0
    $$
    then shows $h^1(X,M\otimes N^{\vee}) \geq h^1(E,\O_E) = q+1$ and $h^2(X,M\otimes N^{\vee}) = h^2(X,M\otimes N^{\vee}\otimes \omega_X) = 0$. Using Riemann--Roch and adjunction, we find
    $$
    h^0(X,M\otimes N^{\vee}) \geq g+1-2d.
    $$
    The non-scalar automorphisms of $F$ have dimension $1 + h^0(X,M\otimes N^{\vee})$. By Lemma~\ref{lem:constancyofLMsections}, $h^0(X,F) = g-d+3$, so
    $$
    \dim \mathcal{W}\leq \dim \mathrm{Gr}(2,H^0(X,F)) -1-h^0(X,M\otimes N^{\vee})\leq g.
    $$
    If $M= N$, then $g=2d+1$ and the non-scalar automorphisms of $F$ have dimension $3$, so we get the same bound $\dim \mathcal{W} \leq g$. Since $g = \dim |\mathcal{L}|$ and the right hand side in the linear growth bound \eqref{eqn:lineargrowthcondition} is always $\geq 0$, this completes the argument.
\end{proof}

\begin{proof}[Proof of Theorem~\ref{thm:intro}]
    Combine Propositions~\ref{prop:genericgreen1}, \ref{prop:genericgreen2} and~\ref{prop:genericgreen3} with the reduction in Proposition~\ref{prop:mainreduction}.
\end{proof}

\noindent It turns out that an anticanonical divisor is not necessary in Theorem~\ref{thm:intro} and it often suffices to assume $-K_X|_C \geq 0$.

\begin{corollary}\label{cor:maincorollary}
    Let $X$ be a birationally ruled surface of irregularity $q$ and let $C \subset X$ be a smooth, non-hyperelliptic curve with $g(C)\geq 3$ and $-K_X|_C \geq 0$. Denote by $c$ the Clifford index of the general curve in $|\mathcal{L}|_{\mathrm{sm}}$. If $q > 0$, denote by $f \colon X \to B$ the ruling and assume $n=\mathrm{deg}(f|_C) \geq 3$. Then, Green's conjecture holds for $C$ if one of the following conditions is satisfied:
    \begin{enumerate}[(i)]
        \item $q = 0$;
        \item $q\geq 1$ and  $n \gon(B)\neq c+2$;
        \item $q \geq 2$, $n  \gon(B)= c+2$ and $B$ has only finitely many gonal pencils.
    \end{enumerate}
\end{corollary}

\begin{proof}
    We show that $X$  admits a birational map $X\to Y$ which restricts to an isomorphism on $C$ and where $-K_Y>0$. The statement then follows directly from Theorem~\ref{thm:intro}. Let $X \to Y$ be a composition of contractions of $(-1)$-curves disjoint from $C$, in some order until no such $(-1)$-curves remain. We claim $-K_Y>0$. Since $Y$ has negative Kodaira dimension, $\omega_Y \not\simeq \O_Y$, so it suffices to show that $\omega_Y^{\vee}$ has a global section. Set $L \coloneqq C+K_Y$ and consider the restriction sequence
    $$
    0\to \O_Y(-L)\to \O_Y(-K_Y)\to \O_C(-K_Y)\to 0.
    $$
    If $L$ is big and nef, then by Kawamata--Viehweg vanishing, restriction induces a surjection 
    $$
    H^0(Y,\O_Y(-K_Y)) \twoheadrightarrow H^0(C,\O_C(-K_Y)).
    $$
    As the contraction map is an isomorphism in a neighborhood of $C$, $H^0(C,\O_C(-K_Y)) \neq 0$ since $-K_X|_C \geq 0$. It follows that $-K_Y$ is effective, so it remains to show that $L$ is big and nef. Identifying $C$ with its image in $Y$, we have $-K_Y.C = -K_X.C \geq 0$, so
    $$
    C^2 = 2g(C)-2-K_Y.C>0,
    $$
    i.e., $C$ is big and nef on $Y$. We first claim that $L$ is nef. Assume for contradiction that there exists an irreducible effective divisor $\Gamma$ with $L.\Gamma <0$. Note that $L$ is effective as $h^0(X,L) \geq  g(C) -q$ from the restriction sequence, and $g(C)-q>0$ from $g(C)\geq 2$ and Riemann--Hurwitz. Then, necessarily, $\Gamma^2<0$ and by nefness of $C$ we have $C.\Gamma \geq 0$. Adjunction then shows $g(\Gamma) = 0$, $\Gamma^2 = -1$ and $\Gamma.K_Y = -1$. It follows 
    $0> L.\Gamma = C.\Gamma -1$ and consequently $C.\Gamma = 0$. Thus $C \cap \Gamma = \emptyset$, a contradiction because $Y$ contains no $(-1)$-curves disjoint from $C$. Assume now that $L$ is not big, i.e., $L^2 = 0$. As $L-K_Y = C$ is big and nef, \cite[Thm. 3.3]{KollarMoriBirationalGeometry} implies that $L$ is semi-ample. Let $\varphi \colon Y\to Z\subset \P^N$ be the map induced by some basepoint-free multiple of $L$. Since $L^2=0$, the image $\varphi(Y) = Z$ is a curve and we assume $Z$ to be smooth (else replace it by its normalization). Let $\varphi'$ be the first map in the Stein factorization of $\varphi$, i.e., $\varphi$ factors as
    \begin{equation*}
        Y \xrightarrow{\varphi'} Z' \rightarrow Z,
    \end{equation*}
    where $Z' \rightarrow Z$ is a finite map between curves and $\varphi'$ has connected fibers. Let $F$ be the numerical class of a fiber of $\varphi$ and $F'$ the numerical class of a fiber of $\varphi'$. Then $F = \deg(Z' \to Z) F'$. Moreover, since $Z$ is a curve, a basepoint-free multiple of $L$ is numerically equivalent to a multiple of $F$. Therefore, $0 = F.L = F.C + F.K_Y$, and the same is true for $F'$ in place of $F$. As $F$ is nef and $C^2 > 0$, we have $F.C > 0$ by the Hodge index theorem. Thus, $F.K_Y < 0$ and the same is true for $F'$. Adjunction now yields that every fiber of $\varphi'$ has arithmetic genus $0$. Hence, $\varphi'$ is a ruling and we get $C.F' = -K_Y.F' = 2$. If $q \geq 1$, then a ruling is unique up to automorphisms of the base curve, so $C.F' = \deg(\varphi'|_C) = \deg(f|_C) \geq 3$, contradicting $C.F' = 2$. If $q = 0$, then $g(Z') = 0$, so $C$ is hyperelliptic. It follows that $L$ is big and nef, as desired.
\end{proof}

\stoptoc 
\resumetoc 

\bibliographystyle{amsalpha} 
\bibliography{TEMPLATEBIB}
\end{document}